\documentclass[a4paper]{article} 
\usepackage[francais]{babel} 
\usepackage[T1]{fontenc} 
\usepackage[latin1]{inputenc}
\usepackage{amsmath}
\usepackage{amssymb}
\usepackage{eufrak}
\usepackage{amsthm}
\usepackage{amsfonts}
\title{Une r\'egion explicite sans z\'eros pour la~fonction~$\zeta$ de Riemann}
\author{Habiba Kadiri} 
\newcommand{\myref}[1]{(\ref{#1})}
\newtheorem{thm}{Th\'eor\`eme}[section]

\newtheorem{lem}[thm]{Lemme}
\newtheorem{prop}[thm]{Proposition}

\begin{document}
\maketitle
{\def\thefootnote{}
}
\begin{abstract}
\noindent
Dans cet article, nous montrons que la fonction $\zeta$ de Riemann n'a
pas de z\'eros dans la r\'egion~:
$$
\Re s \ge 1- \frac1{5.70176 \log |\Im s|} \quad \quad (|\Im s| \ge 2).
$$
\end{abstract}
\section{Historique et r\'esultats.}
\noindent
Depuis l'article de Riemann en 1860 (cf. \cite{Rie}), nous savons que
la r\'epartition des nombres premiers est \'etroitement li\'ee \`a la
r\'epartition des z\'eros d'une fonction particuli\`ere, appel\'ee
depuis la fonction $\zeta$ de Riemann. Nous rappelons qu'elle est
d\'efinie sur le demi-plan $\Re s >1$ par~:
$$
\zeta(s) = \sum_{n\ge1} \frac1{n^s} = \prod_p \left(1-p^{-s}\right)^{-1}
$$
o\`u le produit porte sur les nombres premiers.
\\ \noindent
La s\'erie ci-dessus converge absolument et uniform\'ement dans le
demi-plan $\Re s \ge \sigma_0$, pour tout $\sigma_0 >1$.
La fonction se prolonge en une fonction holomorphe dans le plan
complexe sauf en $1$, qui est pôle unique et en lequel elle a pour
r\'esidu 1. Elle v\'erifie l'\'equation fonctionnelle suivante sur le
plan complexe tout entier~:
$$
\pi^{-s/2} \Gamma(s/2) \zeta(s) = \pi^{-(1-s)/2} \Gamma((1-s)/2) \zeta(1-s)
$$
On en d\'eduit que les z\'eros r\'eels de la fonction $\zeta(s)$ sont les
pôles de $\Gamma(s/2+1)$, c'est \`a dire les entiers $-2n$, o\`u $n\in
\mathbb{N}^*$ et qu'elle a une infinit\'e de racines complexes dont la partie
r\'eelle est comprise entre 0 et 1. Nous appellerons $Z(\zeta)$
l'ensemble de ces z\'eros dits non-triviaux et les noterons $\varrho =
\beta + i \gamma$. Ils se r\'epartissent sym\'etriquement par rapport
\`a l'axe r\'eel et par rapport \`a l'axe $\Re s = 1/2$.
\\ \noindent
L'hypoth\`ese de Riemann affirme qu'en fait, ils se trouvent tous sur
la droite $\Re s = 1/2$. Mais cette conjecture n'a encore \'et\'e ni
d\'emontr\'ee ni contredite.
\\ \noindent
Van de Lune, te Riele et Winter l'ont cependant v\'erifi\'ee en 1986
 (cf. \cite{LRW}) pour les z\'eros de partie imaginaire
 inf\`erieure \`a $5\cdot10^8$, ce qui concerne les $1.5\cdot10^9$
 premiers z\'eros de la fonction $\zeta$ de Riemann. Ce r\'esultat vient même
 d'être tout r\'ecemment am\'elior\'e par S.Wedeniwski jusqu'\`a une partie
 imaginaire de $3\ 330\ 657\ 430.697$.
\\ \noindent
En attendant, l'\'ecriture de $\zeta$ sous forme de produit eul\'erien
nous assure qu'elle n'a pas de z\'eros dans le demi-plan $\Re s >1$. En
fait, l'influence de la formule d'Euler s'\'etend m\^eme \`a gauche de
cette r\'egion. Ainsi, en $1896$, Hadamard (voir \cite{Had}) et De La Vall\'ee Poussin
(cf. \cite{VP1}) \'etablissent simultan\'ement mais s\'epar\'ement que
$\zeta$ ne s'annule pas sur la droite $\Re s = 1$. Cette affirmation
est l'outil fondamental leur permettant d'\'etablir le th\'eor\`eme des nombres premiers, \`a
savoir qu'on a l'estimation asymptotique suivante pour le nombre $\pi(x)$ d'entiers premiers
inf\'erieurs \`a $x$ ~:
$$
\pi(x) \underset{x \rightarrow +\infty}{\sim} Li(x) \ \hbox{ o\`u }\ Li(x) = \int_2^x
\frac{d\,t}{\log t}.
$$
\\ \noindent
En 1899, De La Vall\'ee Poussin (cf. \cite{VP2}) \'elargit son
r\'esultat \`a la r\'egion ~:
$$
\Re s \ge 1- \frac{1}{R_0 \log |\Im s|}, \  |\Im s| \ge 2,
\quad \hbox{avec} \quad R_0=34.82
$$
ce qui lui permet d'estimer le terme d'erreur pour le th\'eor\`eme des
nombres premiers ~:
$$
\pi(x) - Li(x) = \mathcal{O}\left(x \exp \left(-\sqrt{\frac{\log
    x}{R_0}}\right)\right) \quad \hbox{ quand } x \rightarrow +\infty.
$$
\\ \noindent
B.Rosser am\'eliore la valeur de la constante $R_0$, notamment en modifiant le
polynôme trigonom\'etrique (voir ci-apr\`es) et
en $1939$ il obtient $R_0=19$ (cf. \cite{Ros}) puis avec
L.Schoenfeld en 1962, $R_0=17.516$ (cf. \cite{RS1}).
En $1975$, dans \cite{RS2}, ces derniers reprennent une id\'ee
fondamentale d\^ue \`a Stechkin (cf. \cite{Stechkin1}, lemme 2) et
descendent jusqu'\`a $R_0=9.645908801$.
\\ \noindent
Beaucoup plus r\'ecemment, K.Ford
(cf. \cite{Ford}, Th\'eor\`eme 4) atteint une valeur de $8.463$, esentiellement en utilisant une majoration de la fonction $\zeta$ sur
l'axe critique, m\'ethode qui ne se g\'en\'eralise pas aux fonctions
$L$ de Dirichlet.
\bigskip
\\ \noindent
En se basant sur la majoration
de $\zeta(s)$ lorsque $\Re s =1$ donn\'ee par la m\'ethode de
Korobov et Vi\-no\-gra\-dov (cf. \cite{Ko}), il est possible d'obtenir
une r\'egion sans z\'eros du type~:
$$
\Re s > 1- \frac{1}{R_1 (\log |\Im s|)^{2/3} (\log \log |\Im
  s|)^{1/3} } \qquad (|\Im t| \ge 10)
$$
En $1994$, O.V.Popov (cf. \cite{P}) trouve $R_1=14518$,
r\'esultat derni\`erement am\'elior\'e par K.Ford ~:~ $R_1=57.54$ (cf. \cite{Ford}).
\bigskip
\\ \noindent
Nous allons d\'emontrer le r\'esultat suivant~:
\begin{thm}[Principal]\label{THM}\ \\ \noindent
La fonction $\zeta$ de Riemann ne s'annule jamais dans la r\'egion
suivante ~:
$$
\Re s \ge 1-\frac1{R_0 \log \big(|\Im s|\big)},\ |\Im s| \ge 2, \quad
\hbox{ avec }\quad R_0 = 5.70176.
$$
\end{thm}
\ \\ \noindent
Cette r\'egion reste plus large que la r\'egion de
Ford-Vinogradov jusqu'\`a des valeurs de $|\Im s|$ inf\'erieures \`a $e^{9402.562}$.\\ \noindent
D'autre part, les outils mis en \oe uvre ici se g\'en\'eralisent
aux fonctions $L$ de Dirichlet (cf. \cite{kadiri1}).
\\ \\ \noindent
Rappelons les trois points fondamentaux autour desquels s'articule la
preuve d'un tel r\'esultat.
Tout repose d'abord sur une expression de la partie r\'eelle de
$-(\zeta'/\zeta)(s)=\sum_{n\ge1}\Lambda(n)n^{-s}$ en fonction des
z\'eros de la fonction
$\zeta$ de Riemann. Pour cela, il y a deux approches~:
\begin{enumerate}
\item celle, dite globale, de De La Vall\'ee Poussin qui regarde tous
  les z\'eros avec la relation
\begin{equation}\label{1.0}
\Re \Big( \frac{-\zeta'}{\zeta}(s) \Big) =  - \frac12 \log \pi
+ \frac12 \Re \frac{\Gamma'}{\Gamma} \Big(\frac{s}2+1\Big) + \Re
\Big(\frac1{s-1} \Big)
- \sum_{\varrho \in Z(\zeta)} \Re \Big(\frac1{s-\varrho} \Big).
\end{equation}
\item celle, dite locale, de Landau qui s'attache aux z\'eros
  ``proches'' de $s$ avec la majoration
\label{1.0bis}
\begin{eqnarray}
\Re \Big( -\frac{\zeta'}{\zeta}(s) \Big) =  - \sum_{|s-\varrho| \le
  \frac{c}{\log \Im s} } \Re \Big( \frac1{s-\varrho} \Big) +
\mathcal{O} \big( \log \Im s \big)
\end{eqnarray}
\end{enumerate}
Le second point essentiel de la preuve consiste en la positivit\'e de
la somme sur les z\'eros $\sum_{\varrho} \Re \Big(\frac1{s-\varrho}
\Big)$ lorsqu'on suppose $\Re s > 1$.
\\ \noindent
Enfin, la preuve s'ach\`eve avec un autre argument de positivit\'e
~:~\\ \noindent
si $P(\theta) = \sum_{k=0}^K a_k \cos(k\theta)$ est un polyn\^ome
trigonom\'etrique v\'erifiant
\begin{gather*}
a_k \ge 0 \quad \hbox{ et } \quad  P(\theta) \ge 0,\\
\hbox{alors on a ~:}\qquad
\Re \sum_{k=0}^K a_k \sum_{n\ge1} \frac{\Lambda(n)}{n^{\sigma+ikt}} \ge 0.\qquad\qquad\qquad\qquad
\end{gather*}
\\ \noindent
En ce qui concerne le travail pr\'esent\'e ici, nous reprendrons une
id\'ee exploit\'ee entre autres par Stechkin et Heath-Brown, et qui
consiste \`a multiplier la fonction de von Mangoldt $\Lambda(n)$ par une
fonction lisse positive $g(\log n)$ tout en conservant la
positivit\'e de la somme sur les z\'eros. Par exemple
Stechkin (cf. \cite{Stechkin1}) prend
$g(x) = 1-\kappa e^{-\delta x}$ o\`u $\kappa$ et $\delta$ sont deux r\'eels positifs.
Quant \`a Heath-Brown, il propose une formule plus compliqu\'ee
(cf. \cite{HB2} et le paragraphe \ref{Section2}), l'essentiel pour $g$
\'etant d'être de classe $C^2$ dans $]0,+\infty[$, \`a support compact
et, comme nous allons le voir, \`a transform\'ee de Laplace positive.
En fait, nous consid\`ererons le produit de leurs deux fonctions.
\\ \noindent
Les formules de Weil (voir \cite{Weil}) s'appliquent alors et donnent
une formule impliquant la somme sur tous les z\'eros (voir le paragraphe
\ref{Section1})~:
\begin{multline*}
\Re \Bigg( \sum_{n\ge1} \frac{\Lambda(n)}{n^{s}}g(\log n) \Bigg)
=  g(0) \  \Re \Bigg(  - \frac12 \log{\pi} + \frac12
\frac{\Gamma'}{\Gamma} \Big(\frac{s}2+1\Big) \Bigg)
\\  + \Re G(s-1) - \sum_{\varrho \in Z(\zeta)} \Re G(s-\varrho) + \Re R(s)
\end{multline*}
o\`u $G$ est la transform\'ee de Laplace de $g$
$$
G(z) = \int_{0}^{+\infty} e^{-zt} g(t) \,dt = \frac{g(0)}{z} +
\frac1z\int_{0}^{+\infty} e^{-zt} g'(t) \,dt
$$
et o\`u $R(s)$ est un terme reste sous la condition
$G(z) - g(0)/z = \mathcal{O}(1/|z|^2)$.
\\ \noindent
En remarquant que la transform\'ee de Laplace de la fonction
constante \'egale \`a 1 est $\frac1z$, et en rappelant que
l'in\'egalit\'e $\Re \frac1z \ge 0$ si $\Re z > 0 $ est fondamentale
pour traiter la somme sur les z\'eros, nous souhaiterions imposer \`a
$g$ que sa transform\'ee de Laplace v\'erifie~:
$$
\Re{G(z)} \ge 0 \quad \hbox{si} \quad \Re z > 0
$$
\noindent
Nous affaiblissons ici cette condition  en utilisant la sym\'etrie des
z\'eros non triviaux et
g\'en\'eralisons le lemme de Stechkin (cf. lemme 2, \cite{Stechkin1}) en montrant que
$$
\Re G(s-\varrho) + \Re G(s+1-\overline \varrho) \ge 0\qquad \hbox{si}\quad \Re(s-\varrho)>0
$$
(voir le paragraphe \ref{Section4}).
\\ \noindent
En prenant $\Re s > 1$, une telle d\'emonstration apporterait d\'ej\`a
une am\'elioration \`a la constante de Rosser et Schoenfeld. Or la
fonction $g$ \'etant \`a support compact, cela nous permet de choisir
$s$ avec $\Re s \le 1$.
\\ \noindent
Pour ce qui est de la somme
$$
\sum_{{\begin{substack}{\varrho \in Z(\zeta) \\
                       \Re(s-\varrho)\le0}
      \end{substack}}}
 \Re G(s-\varrho),
$$
nous montrerons que c'est un terme reste (voir le paragraphe \ref{Section4}).
\\ \noindent
Enfin, nous conservons l'argument trigonom\'etrique final en
consid\`erant un polynôme de degr\'e 4 proche de celui introduit par Rosser et
Schoenfeld (voir le paragraphe \ref{Section3}).
\bigskip
\\ \noindent
Nous donnons dans le paragraphe qui suit tous les r\'esultats
n\'ecessaires \`a l'\'etablissement de notre r\'esultat et nous y fixons nos
notations. Pour le
d\'etail des preuves, nous nous reporterons ensuite aux parties trois
et quatre.
\bigskip
\\ \noindent
Je remercie O. Ramar\'e pour les conseils avis\'es
qu'il m'a prodigu\'es au fil de cette \'etude ainsi que K. Ford pour
m'avoir aimablement transmis une version pr\'eliminaire de son article
\cite{Ford}.
\\ \noindent
Enfin, je tiens \'egalement \`a remercier le r\'ef\'er\'e pour tout le soin qu'il a
apport\'e \`a la relecture du manuscrit. Je lui suis particuli\`erement
reconnaissante d'avoir trouv\'e une optimisation au polynôme de Rosser et
Schoenfeld et au param\`etre $\theta$, ce qui permet d'am\'eliorer sensiblement
le r\'esultat final.

\section{Structure de la preuve.}
\noindent Commen\c cons par pr\'eciser nos param\`etres.
Nous nous donnons
une fonction positive $f$, de classe $C^2([0,d])$, à support
compact dans $[0,d[$ et telle que :
\begin{gather}
f(d)=f'(0)=f'(d)=f"(d)=0.
\tag{$H_1$}
\end{gather}
dont nous notons $F$ la transformée de Laplace~:
$$
F(s)=\int_0^d{e^{-st}f(t)}\,dt.
$$
Cette fonction sera choisie au paragraphe \ref{Section2}.
\noindent
Nous consid\'erons ensuite un z\'ero non trivial $\varrho_0 =
\beta_0+i\gamma_0$ de la
fonction $\zeta$ de Riemann et nous souhaitons montrer que ce
z\'ero v\'erifie le th\'eor\`eme \ref{THM}.
La sym\'etrie des z\'eros de $\zeta$ nous permet de nous limiter au cas
o\`u $\gamma_0$ est positif.
De plus, comme tous les z\'eros de partie imaginaire positive inf\'erieure
\`a $T_0 = 3.3 \cdot 10^9$ sont connus (cf. \cite{Wed}) et r\'esident tous sur la
droite critique, nous supposerons que
$\gamma_0$ est sup\'erieur \`a $T_0$. Enfin, nous supposerons que
$(1-\beta_0)\log\gamma_0 \le \frac15$
\noindent
Nous noterons $\eta$ le r\'eel $(1-\beta_0)$, $s$ le
nombre complexe $\sigma+it$, o\`u $t \in [0,+\infty[$, $R$ un r\'eel
pour lequel la r\'egion sans z\'ero est v\'erifi\'ee et $t_0$ un r\'eel sup\`erieur \`a $1$.
Nous \'ecrirons $\eta $ sous la forme
$\frac1{r\log\gamma_0}$, o\`u $5 \le r \le R$ en vertu de notre hypoth\`ese, et $\sigma$ sous la forme $1-\frac1{R\log(4\gamma_0+t_0)}$.  Grâce au
r\'esultat de Rosser (cf. \cite{Ros}), nous prenons tout d'abord $R=9.645908801$.
Notamment, $\sigma \ge \sigma_0=1-\frac{1}{9.645908801\log(4T_0+1)}\ge0.99555$ et $ \eta
\le \eta_0 = \frac1{r\log T_0} \le \frac1{5\log T_0}\le0.00913$.
\\ \noindent
Dans la suite, $\kappa$ et $\delta$ d\'esignent des
constantes qui d\'ependent et ne d\'ependent que de $r$ et $R$. Cette
d\'ependance est assez faible mais toutefois num\'eriquement
interessante. De plus, nous leur imposons la condition suivante ~:
$$
\Big(\delta^{-3}+\big(1-\eta_0+\delta\big)^{-3}\Big)^{-1} \le \kappa
\le \Big(\delta^{-1}+\big(1-\eta_0+\delta\big)^{-1}\Big)^{-1}.
$$
Ou plut\^ot, en fixant les valeurs de $\eta_0$ dans $[0;10^{-2}]$ et
$\delta$ dans $[(\sqrt5-1)/2;0.866]$, nous demandons \`a $\kappa$ de
v\'erifier ~:
\begin{equation}
\label{hyp-kappa}
\Big(\delta^{-3}+\big(1+\delta\big)^{-3}\Big)^{-1} \le \kappa
\le \Big(\delta^{-1}+\big(0.99+\delta\big)^{-1}\Big)^{-1}.
\end{equation}
\subsection{Une formule explicite.}
\label{Section1}\noindent
Nous commen\c cons par une formule explicite \`a la Weil
(cf. \cite{Weil})
que nous d\'emontrons au paragraphe \ref{Section5.1}.
\begin{prop}
\label{Prop1.1}
Soit \(f\) une fonction comme ci-dessus et soit $s$ un nombre complexe. Nous avons
\begin{multline}
\label{1.1}
\Re \Bigg(\sum_{n\ge1} \frac{\Lambda(n)}{n^{s}} f(\log n) \Bigg) =
 f(0) \Bigg( - \frac12 \log \pi +\Re \frac12 \frac{\Gamma'}{\Gamma}\Big(\frac{s}2+1\Big) \Bigg)
\\+\Re F(s-1) -\sum_{\varrho\in Z(\zeta)}\Re F(s-\varrho)
\\ + \Re \Bigg( \frac{1}{2i\pi}\int_{1/2-i\infty}^{1/2+i\infty}
\Re\frac{\Gamma'}{\Gamma}\Big(\frac z2 \Big)
\frac{F_2(s-z)}{(s-z)^2}d\,z + \frac{F_2(s)}{s^2}
\Bigg)
\end{multline}
o\`u $F_2$ est la transform\'ee de Laplace de $f"$ et o\`u $Z(\zeta)$
d\'esigne l'ensemble des z\'eros non triviaux de $\zeta$.
\end{prop}\noindent
Nous prenons des notations suppl\'ementaires pour all\'eger
quelque peu le travail typographique et posons
\begin{eqnarray*}
T_1(s) &=& - \frac12 \log(\pi)
+ \frac12 \Re \frac{\Gamma'}{\Gamma}\Big(\frac{s}2+1\Big) ,
\\T_2(s) &=&
\Re\Bigg(\frac{1}{2i\pi}\int_{1/2-i\infty}^{1/2+i\infty}
\Re\frac{\Gamma'}{\Gamma}\Big(\frac z2 \Big)
\frac{F_2(s-z)}{(s-z)^2}d\,z + \frac{F_2(s)}{s^2} \Bigg)
\nonumber\\ &=&
\frac{1}{2\pi}\int_{-\infty}^{+\infty}
\Re\frac{\Gamma'}{\Gamma}\Big(\frac14 + i\frac t2 \Big) \Re
\frac{F_2(s-1/2-it)}{(s-1/2-it)^2}d\,t + \Re \frac{F_2(s)}{s^2}.
\end{eqnarray*}
Nous introduisons aussi les trois diff\'erences
\begin{gather*}
\Delta_{1}(s) = T_1(s)\ -\ \kappa\ T_1(s+\delta)
\quad,\quad
\Delta_{2}(s) = T_2(s)\ -\ \kappa\ T_2(s+\delta)
\\
D(s)=\Re F( s)-\ \kappa\ \Re F( s+\delta).
\end{gather*}
\noindent
Notons qu'il est sous-entendu qu'elles d\'ependent des param\`etres
$\delta$ et $\kappa$.
\\ \\ \noindent
De \myref{1.1}, nous tirons~:
\begin{multline}
\label{1.2}
\Re \sum_{n\ge1} \frac{\Lambda(n)}{n^{s}} f(\log n)
\Big(1-\frac{\kappa}{n^\delta}\Big)
\\= f(0)\Delta_{1}(s) +
D(s-1) - \sum_{\varrho\in Z(\zeta)}D(s- \varrho)\ +
\Delta_{2}(s).
\end{multline}
\subsection{La fonction test f.}
\label{Section2}\noindent
Nous noterons $\tilde F (X,Y)$ la partie r\'eelle de la transform\'ee
de Laplace de $f$~:
\begin{equation}
\label{DefFtilde}
\tilde F (X,Y) = \Re \int _0^d e^{-(X+iY)t}f (t) d\,t.
\end{equation}
En sus des conditions requises \`a $f$ dans l'introduction, nous
imposons \`a $\tilde F$ de v\'erifier~:
\begin{equation}
\tilde F(X,Y) \ge 0 \quad \hbox{si} \quad X \ge 0.
\tag{$H_2$}
\end{equation}
Heath-Brown propose une famille de fonctions
(cf. lemme~7.5, \cite{HB2})
dont nous ne savons pas
si elles sont optimales pour notre probl\`eme, mais dont nous pensons
qu'elles sont assez bien adapt\'ees \`a notre m\'ethode.
Pour $\theta \in ]\pi/ 2,\pi[$, nous d\'efinissons~:
$$
f(t) = \eta h_\theta(\eta t)
$$
o\`u $h_\theta$ est ind\'ependante de $\eta$. Cette fonction est nulle
en dehors de $[0,{-2\theta}/ { \tan \theta}]$ et, pour $u$ appartenant
\`a cet intervalle, vaut
\begin{multline*}
h_\theta(u) =
(1+\tan^2 \theta)
\Big[ (1+\tan^2 \theta)
\Big(\frac{-\theta}{\tan \theta} - \frac u2 \Big)
\cos(u \tan \theta)+ \frac{-2\theta}{\tan \theta} - u \\
-\frac{\sin(2\theta+ u\tan \theta)}{\sin(2\theta)}
+ 2 \Big( 1+\frac{\sin(\theta+ u \tan\theta)}{\sin\theta}\Big)\Big].
\end{multline*}
Nous avons alors
$$
f(0) = \eta g_1(\theta)
,\
\Re F(0) = \tilde F(0,0) = g_2(\theta)
,\
\Re F(1 - \beta_0) = \tilde F(1-\beta_0, 0) = g_3(\theta)
$$
o\`u nous avons pos\'e~:
\begin{equation*}
\left\{
\begin{array}{rcl}
g_1(\theta)&=&(1+\tan^2\theta)(3-\theta \tan \theta -3 \theta \cot\theta),\\
g_2(\theta)&=&2(1+\tan^2\theta)(1-\theta \cot \theta)^2,\\
g_3(\theta)&=&2\tan^2 \theta + 3 -3 \theta \tan \theta -3 \theta\cot \theta.
\end{array}
\right.
\end{equation*}
Nous prendrons
$\theta=1.848$
ce qui nous donnera
\begin{gather*}
g_1(\theta) = 147.84112+\mathcal{O^*}(10^{-5})
,\quad g_2(\theta) = 62.17067+\mathcal{O^*}(10^{-5}),
\\g_3(\theta) = 48.76676+\mathcal{O^*}(10^{-5})
\end{gather*}
o\`u $u=\mathcal{O^*}(v)$ signifie $|u|\le v$.
\noindent
Pour les besoins ult\'erieurs, nous d\'efinissons aussi
$$
d(\theta,\eta) = \frac{d_1(\theta)}{\eta}\quad \text{o\`u}\quad
d_1(\theta) =\frac{-2\theta}{ \tan \theta}=1.05161+\mathcal{O^*}(10^{-5}).
$$
\subsection{Une in\'egalit\'e trigonom\'etrique.}
\label{Section3}\noindent
Nous utilisons ici l'in\'egalit\'e suivante~:
$$
\sum_{k=0}^4 a_k \cos(ky) = 8(0.91+\cos y)^2(0.265+\cos y )^2 \ge 0,
$$
avec
\begin{gather*}
a_0 = 10.91692658 ,
\quad  a_1 = 18.63362,
\quad  a_2 = 11.4517,
\\ a_3 = 4.7,
\quad a_4 = 1,
\quad A = \sum_{k=1}^4 a_k = 35.78532.
\end{gather*}
En remarquant que
$$
 \sum_{n\ge1}f ( \log n ) \frac{\Lambda(n)}{n^{\sigma}}
\Big(1-\frac{\kappa}{n^\delta}\Big) \
 \sum_{k=0}^4 a_k \cos(k\gamma_0 \log n) \ge 0
$$
nous obtenons grâce \`a \myref{1.2} l'in\'egalit\'e fondamentale suivante~:
\begin{multline}
\label{1.3}
\sum_{k=0}^4 a_k \Bigg( f(0)\Delta_{1}(\sigma+ik\gamma_0) + D(\sigma-1+ik\gamma_0)
\\ - \sum_{\varrho\in Z(\zeta)} D(\sigma+ik\gamma_0- \varrho)\
+ \Delta_{2}(\sigma+ik\gamma_0) \Bigg) \ge 0
\end{multline}
Il reste donc \`a trouver des majorations pour chacun des termes
ci-dessus. C'est l'objet du paragraphe suivant.
\subsection{La r\'egion sans z\'eros.}
\label{Section4}\noindent
Les valeurs que nous donnons ici sont calcul\'ees pour $r=5.97484$, $R=9.645908801$. A
la fin de ce paragraphe, la valeur $R=5.97485$ est donc licite et nous pouvons
recommencer les calculs. Ce que nous avons fait, mais l'\'etape principale est
la premi\`ere et elle permet en outre au lecteur de v\'erifier nos r\'esultats.
\noindent
Pour $T_0$, nous prendrons la valeur de S.Wedeniwski, c'est \`a dire tr\`es
exactement $3\, 330\, 657\, 430.697$. (La valeur finale obtenue pour $R_0$ sera alors
am\'elior\'ee d'un centi\`eme par rapport \`a la valeur $T_0$ de te Riele et
Winter.)

\par
\noindent
Commen\c cons par le terme $\Delta_{1}(s)$~: lorsque $t$ est non nul, $\Re\frac{\Gamma'}{\Gamma}\big(\frac{s}2+1\big)$ est de l'ordre de
$\log t$ (voir \myref{Gamma3} au paragraphe \ref{Section5.3}). Au paragraphe \ref{Section2.3} nous \'etablirons la proposition suivante~:
\begin{prop}
\label{Prop1.2}
Il existe une fonction $\mathcal{C}_1(\eta)$ qui v\'erifie
$$
f(0) \sum_{k=0}^4 a_k \Delta_{1} (\sigma + ik\gamma_0) \le \frac A2
(1-\kappa) g_1(\theta) \eta \log \gamma_0 + \mathcal{C}_1(\eta).
$$
On se reportera \`a
\myref{C1} et au lemme \ref{delta1}, pour la d\'efinition de $\mathcal{C}_1(\eta)$ et
$$
\mathcal{C}_1(\eta)\le - 2718.913\, \eta.
$$
\end{prop}
\par\noindent
En ce qui concerne le terme $D(s)$,
nous avons $\Re F(s-1)=\tilde F(\sigma-1,0)$ lorsque $t=0$.
Sinon $\Re F(s-1)$ est de l'ordre de $\eta/t^2$ (voir la proposition \ref{Prop2.5}) et nous montrerons au paragraphe \ref{Section2.4} que
\begin{prop}
\label{Prop1.3}
Il existe une fonction $\mathcal{C}_2(\eta)$ qui v\'erifie
$$
\sum_{k=0}^4a_k D(\sigma-1+ik\gamma_0) \le a_0 \tilde F(\sigma-1,0) +
\mathcal{C}_2(\eta)
$$
La fonction $C_2(\eta)$ est d\'efinie en
\myref{C2} et
$$
\mathcal{C}_2(\eta)\le -1\,141.389 \,\eta + 2.794 \cdot10^{-15} \,\eta^2 + 26\,515.117 \,\eta^3.
$$
\end{prop}
\par\noindent
Le terme $\Delta_{2}(s)$ est un terme reste, de l'ordre de
$\eta^3$ (voir le lemme \ref{LemmeD2}). Nous montrerons
au paragraphe \ref{Section2.6} la proposition suivante~:
\begin{prop}
\label{Prop1.4}
Il existe une fonction $\mathcal{C}_4(\eta)$ qui v\'erifie
$$
\sum_{k=0}^4 a_k \Delta_{2}(\sigma+ik\gamma_0) \le
\mathcal{C}_4(\eta)
$$
Les \'el\'ements qui d\'efinissent $\mathcal{C}_4(\eta)$ sont donn\'es en \myref{C4}, \myref{C41} et \myref{C42} et
$$
\mathcal{C}_4(\eta)\le 2.3887 \cdot10^6\,\eta^3
$$
\end{prop}
\par\noindent
Nous en arrivons enfin au point essentiel, trait\'e au paragraphe \ref{Section2.5}, qui est l'\'etude de la
somme sur les z\'eros.
\noindent
Pour cela, en nous inspirant de l'id\'ee de Stechkin bas\'ee sur la sym\'etrie
des z\'eros de $\zeta$, nous pouvons r\'e\'ecrire la premi\`ere somme
sous la forme~:
\begin{multline*}
\sum_{\varrho\in Z(\zeta)}D(s-\varrho)
=  D(s-\varrho_0) + D(s-1+\overline{\varrho_0})
\\ + \frac12 \sum_{\varrho\in Z(\zeta) \setminus
  \{\varrho_0,1-\overline\varrho_0\} } \Big[
D(s-\varrho)+D(s-1+\overline{\varrho}) \Big]
\end{multline*}
La  proposition \ref{Prop2.6} nous permettra d'\'eliminer une partie des termes de la somme grâce au r\'esultat suivant~:
\begin{gather*}
 D(s-\varrho) + D(s-1+\overline{\varrho}) \ge 0
\\ \text{\ si\ } 1-\sigma < \beta < \sigma  \text{\ ,\ \ }
 \kappa = \kappa_0 \text{\ \ et\ \ }
 \delta = \delta_0
\end{gather*}
o\`u $\delta_0$ est la solution de l'\'equation $\kappa_2(\delta)=\kappa_3(\delta)$, $\kappa_2(\delta)$ et $\kappa_3(\delta)$ \'etant respectivement d\'efinis en \myref{k2} et \myref{k3}, et o\`u $\kappa_0$ est la valeur de $\kappa_2$ en $\delta_0$.
\noindent
\`A un $\mathcal{O}(\eta_0)$ pr\`es, nous avons en fait que
\begin{gather*}
\kappa_2(\delta) = \frac1{1+2\delta} \quad \hbox{et} \quad \kappa_3(\delta) = \frac1{\frac1{\delta}+\frac1{1+\delta}}
\end{gather*}
ce qui permet d'approcher $\delta_0$ par $\frac{\sqrt5-1}2$ et $\kappa_0$ par $\frac1{\sqrt5}$.
\noindent
Plus exactement, nous trouvons $\delta_0=0.62063+\mathcal{O}^*(10^{-5})$ et
$\kappa_0=0.4389+\mathcal{O}^*(10^{-5})$ pour $r=5.97484$.
\noindent
Il reste alors \`a minorer la somme restante qui porte sur les z\'eros de partie r\'eelle v\'erifiant $\sigma\le\beta\le1$, ce qui revient \`a $\gamma\ge \gamma_0+t_0$. Cette somme se trouve alors d\'ependre de la valeur $t_0$. Nous verrons \`a la proposition \ref{Prop2.7} comment la minorer et nous obtiendrons finalement~:
\begin{prop}
\label{Prop1.5}
Il existe une fonction $\mathcal{C}_3$ qui v\'erifie
$$
\sum_{k=0}^4 a_k\sum_{\varrho \in Z(\zeta)} D(\sigma +ik
\gamma_0 - \varrho) \ge a_1 \tilde F(\sigma-\beta_0,0) - \mathcal{C}_3(\eta)
$$
Pour une d\'efinition explicite de $\mathcal{C}_3$, nous nous
reporterons \`a \myref{C3bis}. En attendant, nous pouvons toujours voir que~:
$$
\mathcal{C}_3(\eta) \le   -54.957\, \eta +
344\,602.065\, \eta^2 + 3\,384\,045.191\,\eta^3.
$$
\end{prop}
\par\noindent
Finalement, en notant $\mathcal{C}(\eta) =
\mathcal{C}_1(\eta) + \mathcal{C}_2(\eta) + \mathcal{C}_3(\eta) +
\mathcal{C}_4(\eta)$, nous tirons de l'in\'egalit\'e fondamentale \myref{1.3}~:
$$
0 \le  \frac{A}{2}(1-\kappa)g_1(\theta) \eta \log \gamma_0
 + a_0 \tilde F(\sigma-1,0)
-  a_1 \tilde F(\sigma-\beta_0,0)
+ \mathcal{C}(\eta)
$$
soit encore
$$
\eta \log \gamma_0 \ge
\frac{
a_1 \tilde F(\sigma-\beta_0,0) - a_0 \tilde F(\sigma-1,0)
- \mathcal{C}(\eta)}
{\frac{A}2 g_1(\theta)(1-\kappa)}
$$
En fait, $\mathcal{C}(\eta)$ s'\'ecrit sous la forme~:
\begin{gather}
 \alpha_1\, \eta + \alpha_2\, \eta^2 + \alpha_3\, \eta^3,\\
\hbox{ o\`u }\ \alpha_1 = -3915.260\, , \ \alpha_2 = 344\,602.439\, , \ \alpha_3 = 5\,799\,250\,.773.
\end{gather}
Comme $\alpha_1$ est une constante n\'egative et $\alpha_2$ et
$\alpha_3$ deux constantes positives, on voit facilement que
$\mathcal{C}(\eta)$ admet trois racines r\'eelles dont une \'egale
\`a z\'ero et les deux autres de signes oppos\'es.
$\mathcal{C}(\eta)$ est donc successivement n\'egatif puis positif sur
$[0,+\infty[$ et on obtient que $\mathcal{C}(\eta)$ est n\'egatif sur
$[0,\eta_0]$ en v\'erifiant que $\mathcal{C}(\eta_0)$ l'est~:~nous
trouvons $\mathcal{C}(\eta_0)=-7.22827$.
\\ \noindent
La constante cherch\'ee est ainsi donn\'ee par
\begin{equation}
\label{Rfinal}
\frac{\frac{A}2 g_1(\theta)(1-\kappa)}
{a_1 \tilde F(\sigma-\beta_0,0) - a_0 \tilde F(\sigma-1,0)}
\end{equation}
qu'on optimise en $\sigma$.
\noindent
En notant $\omega = \frac{1-\sigma}{\eta}$, le terme $a_1 \tilde F(\sigma-\beta,0) - a_0 \tilde F(\sigma-1,0)$
 s'av\`ere être une fonction de $\omega$ que nous noterons $K$~:
$$
K(\omega) = \int_{0}^{d_1(\theta)} ( a_1 e^{-t} - a_0) h_{\theta}(t) e^{\omega t} d\,t
$$
$K$ est une fonction croissante sur $[0,1]$, sa valeur optimale est donc en $\omega=\frac{r \log T_0}{R
\log(4T_0+1)}$, en vertu des hypoth\`eses faites sur $\sigma$ et $\eta$.
\noindent
Les donn\'ees initiales $r=5.97484$, $R=9.645908801$ nous am\`ene \`a
$R_0=5.97485$. Nous allons optimiser notre r\'esultat en r\'eit\'erant
les calculs~:
rempla\c cons $R$ par la valeur $R_0$ que nous venons de trouver et $r$
par une valeur sup\'erieure \`a celle que nous venons d'utiliser mais
inf\'erieure au futur $R_0$ (nous proc\'edons par tatonnement). Les
valeurs successives de $r$ et $R$ que nous obtenons ainsi forment deux
suites d\'ecroissantes qui semblent tendre vers une valeur commune.
Nous avons choisi de nous arr\^eter \`a une pr\'ecision de
$10^{-5}$ pour la constante $R_0$, c'est \`a dire \`a la
sixi\`eme \'etape. \\ \noindent
Nous donnons ci-dessous les valeurs successives prises pour $r$ et
$R$, ainsi que celles des
param\`etres $\eta_0$, $\kappa$ et $\delta$ impliqu\'es, et enfin celles trouv\'ees pour
$R_0$.
$$
\begin{array}{|c|c|c|c|c|c|c|c|c}
\hline
R & r&\eta_0 \cdot{10^3}&\kappa&\delta
&R_0\\
\hline
9.645908801&5.97484&7.63319&0.438904&0.620626
&5.974849075\\
5.974849075&5.73045&7.95873&0.438525&0.620748
&5.730454010\\
5.730454010&5.70487&7.99441&0.438483&0.620762
&5.704872616\\
5.704872616&5.70208&7.99832&0.438479&0.620763
&5.702089881\\
5.702089881&5.70178&7.99874&0.438478&0.620763
&5.701785245\\
5.701785245&5.70174&7.99880&0.438478&0.620763
&5.701752890\\\hline
\end{array}$$
%
Nous pouvons ainsi prendre $R_0 = 5.70175$. Nous remarquerons que cette valeur est
assez proche de la valeur optimale calcul\'ee en $\omega=\frac{r}{R}$ et qui
vaut $5.65267$.
 \\\\
  \noindent
Les calculs ont \'et\'e men\'es \`a la fois sous MAPLE et sous PARI / GP avec une
pr\'ecision de $10^{-28}$ et dans chacun des cas, nous retrouvons les
r\'esultats annonc\'es. \newline
Nous pr\'ecisons qu' en utilisant le polynôme de Rosser et Schoenfeld, c'est \`a
dire $8(0.9126+\cos y)^2(0.2766+\cos y)^2$, et en prenant $\theta=1.848$, nous
trouvons pour $R_0$ la valeur $5.70216$.\newline
D'autre part, en ce qui concerne le choix de la valeur de $\theta$, nous
remarquons que le terme final \'etudi\'e \myref{Rfinal}
$\frac{A g_1(\theta)}{K(\omega)}$, est en fait une fonction d\'ependant
uniquement des trois param\`etres $r$, $R$ et $\theta$.
Pour chaque \'etape d\'ecrite pr\'ec\'edemment, c'est \`a dire pour
chaque $r$ et $R$ choisi, nous pouvons donc calculer
la valeur de $\theta$ en laquelle \myref{Rfinal} est optimal.\newline
En prenant pour données initiales $r=5.97145$, $R=9.645908801$, nous trouvons
ainsi qu'en $\theta=1.85362+\mathcal{O}^*(10^{-5})$, la valeur $R_0$ vaut
$5.97145+\mathcal{O}^*(10^{-5})$ et en r\'eit\'erant le proc\'ed\'e :
$$
\begin{array}{|c|ccccccc|}
\hline
R & 9.645908801 & 5.97146 & 5.73009 & 5.70484 & 5.70210 & 5.70180 & 5.70176 \\
\hline
r & 5.97145 & 5.73008 & 5.70483 & 5.70208 & 5.70178 & 5.70174 & 5.70174 \\
\hline
\theta & 1.85362 & 1.84834 & 1.84781 & 1.84775 & 1.84774 & 1.84774 & 1.84774 \\
\hline
R_0 & 5.97146 & 5.73009 & 5.70484 & 5.70210 & 5.70180 & 5.70176 & 5.70175 \\
\hline
\end{array}
$$\noindent
Finalement, nous avons  choisi par souci de clart\'e de fixer la valeur de
$\theta$ \`a $1.848$, d'autant plus que cela n'influe pas sur la pr\'ecision donn\'ee au r\'esultat final.
\section{Pr\'eliminaires.}
\label{Section5}\noindent
Cette partie se d\'ecompose elle-m\^eme en deux. Tout d'abord nous
\'etablissons une formule explicite assez g\'en\'erale. Ensuite, nous
\'etudions plus en d\'etails la fonction $\tilde F$ introduite en
\myref{DefFtilde}.
\subsection{Formule explicite.}
\label{Section5.1}
\begin{thm}
\label{FormuleExplicite}
Soit $\phi$ une fonction \`a valeurs complexes d\'efinie sur la droite
r\'eelle qui v\'erifie les conditions (A) et (B) suivantes~:
\begin{enumerate}
\item[(A)] $\phi$ est continue et continuement d\'erivable sur $\mathbb{R}$ sauf en
un nombre fini de points $a_i$ o\`u $\phi(x)$ et sa d\'eriv\'ee
$\phi'(x)$ n'ont que des discontinuit\'es de premi\`ere esp\`ece et
pour lesquels $\phi$ v\'erifie la condition de la moyenne (i.e
$\phi(a_i) = \frac12 [\phi(a_i+0)+\phi(a_i-0)]$).
\item[(B)] Il existe $b>0$ tel que $\phi(x)e^{x/2}$ et $\phi'(x)e^{x/2}$ soient
$\mathcal{O}(e^{-(1/2+b)|x|})$ au voisinage de l'infini.
\end{enumerate}\noindent
Pour tout r\'eel $a<1$, v\'erifiant $0<a<b$, $\phi(x)$ poss\`ede alors
une transform\'ee de Laplace
$$
\Phi (s) = \int_{0}^{+\infty}\phi (x)e^{-sx}d\,x
$$
qui est holomorphe dans la bande $-(1+a)<\sigma<a$ et qui est
$\mathcal{O}(1/|t|)$ uniform\'ement dans la bande $-(1+a) \le \sigma
\le a$.\\\\ \ \noindent
Soient $q$ un entier non nul et $\chi$ un caract\`ere primitif de
Dirichlet modulo $q$.\\ \ \noindent
Notons $\delta_{q,1} = \begin{cases}0 & \mathrm{si\ } q=1 \\ 1
& \mathrm{sinon}\end{cases}$
et $\mathfrak{a}=\begin{cases}0 & \mathrm{si\ } \chi(-1)=1 \\ 1
& \mathrm{sinon}\end{cases}$. Nous avons alors
\begin{multline*}
\sum_{n\ge1}\Lambda(n)\chi(n)\phi(\log n) = \delta_{q,1}
\Big( \Phi(-1) + \Phi(0) \Big) +\frac12 (1-\delta_{q,1})(1-\mathfrak{a})\Phi(0)
\\- \sum_{\varrho\in Z(\chi)}\Phi(-\varrho)+  \phi (0)\log \frac q\pi
+ \sum_{n\ge1}\frac{\Lambda(n)\overline\chi(n)}{n}\phi (-\log n)
\\+ \lim_{T\rightarrow+\infty}\frac1{2i\pi}
\int_{1/2-iT}^{1/2+iT}
\Re\frac{\Gamma'}{\Gamma}\Big(\frac{s+\mathfrak{a}}{2}\Big)\Phi(-s)d\,s
\end{multline*}
\end{thm}
\begin{proof}
\ \newline\noindent
D'apr\`es le th\'eor\`eme d'inversion de Laplace~:
$$
\phi(\log n) = \frac1{2i\pi}
\int_{-(1+a)-i\infty}^{-(1+a)+i\infty}\Phi(s) n^s d\,s \quad(n\ge1).
$$
Ainsi, par d\'efinition de $\frac{L'}{L}(s,\chi)$ pour $\sigma>1$ et grâce
au changement de variable $s\mapsto -s$, nous pouvons \'ecrire~:
\begin{eqnarray}
\label{3.1}
\sum_{n\ge1}\Lambda(n)\chi(n)\phi(\log n)
=\frac1{2i\pi} \int_{1+a-i\infty}^{1+a+i\infty} -\frac{L'}{L}(s,\chi) \Phi(-s) d\,s.
\end{eqnarray}
Soit $T>0$, notons $I(T)$ l'int\'egrale~:
\[ \frac1{2i\pi} \int_{1+a-iT}^{1+a+iT} -\frac{L'}{L}(s,\chi)
\Phi(-s) d\,s.
\]
L'int\'egration de $-\frac{L'}{L}(s,\chi) \Phi(-s)$ sur le contour du
rectangle form\'e par les droites $\sigma=1+a,\sigma=-a,t=T,t=-T$
permet de r\'e\'ecrire $I(T)$~:
\begin{multline}
\label{3.2}
I(T) = \frac1{2i\pi} \int_{-a-iT}^{-a+iT}
-\frac{L'}{L}(s,\chi) \Phi(-s) d\,s
\\+ \frac1{2i\pi} \int_{-a+iT}^{1+a+iT} -\frac{L'}{L}(s,\chi) \Phi(-s) d\,s
\\ - \frac1{2i\pi} \int_{-a-iT}^{1+a-iT} -\frac{L'}{L}(s,\chi)
\Phi(-s) d\,s
\\ - \Bigg( - \delta_{q,1}\Phi(-1) + \frac12
(1-\delta_{q,1})(1-\mathfrak{a})\Phi(0) + \sum_{\varrho\in
Z(\zeta)}\Phi(-\varrho) \Bigg)
\end{multline}
Gr\^ace \`a la condition (B), les deux derni\`eres int\'egrales
tendent vers $0$ lorsque $T$ tend vers $\infty$.
\noindent
De plus, l'\'equation fonctionnelle de $L$~:
$$
 -\frac{L'}{L}(s,\chi) = \log \frac q\pi - \frac{L'}{L}(1-s,\overline\chi) + \frac12
\Big\{ \frac{\Gamma'}{\Gamma}\Big(\frac{s+\mathfrak{a}}{2}\Big)
+\frac{\Gamma'}{\Gamma}\Big(\frac{1-s+\mathfrak{a}}{2}\Big) \Big\}
$$
permet de d\'ecomposer la premi\`ere int\'egrale en la somme des trois
int\'egrales suivantes~:
\begin{eqnarray*}
I_1(T) &=& \frac1{2i\pi} \int_{-a-iT}^{-a+iT} \log \frac
q\pi \Phi(-s) d\,s
\\I_2(T) &=& \frac1{2i\pi} \int_{-a-iT}^{-a+iT} -\frac{L'}{L}(1-s,\overline\chi) \Phi(-s) d\,s
\\I_3(T) &=& \frac1{2i\pi} \int_{-a-iT}^{-a+iT} \frac12 \Big\{
\frac{\Gamma'}{\Gamma}\Big(\frac{s+\mathfrak{a}}{2}\Big)
+\frac{\Gamma'}{\Gamma}\Big(\frac{1-s+\mathfrak{a}}{2}\Big) \Big\} \Phi(-s) d\,s
\end{eqnarray*}
D'une part, le th\'eor\`eme d'inversion de Laplace permet d'\'ecrire
imm\'ediatement~:
\begin{eqnarray}
\label{3.3}
 \lim_{T\rightarrow +\infty} I_1(T) =  \phi(0)\log \frac q\pi
\end{eqnarray}
D'autre part, le d\'eveloppement de $-\frac{L'}{L}$ en s\'erie de
Dirichlet donne~:
\begin{eqnarray}
\label{3.4}
 \lim_{T\rightarrow +\infty} I_2(T) =
\sum_{n\ge1}\frac{\Lambda(n)\overline\chi(n)}{n}\phi(-\log n)
\end{eqnarray}
En ce qui concerne $I_3$, d\'epla\c cons la droite d'int\'egration vers
la droite
$\sigma=1/2$ sur laquelle $\Gamma$ v\'erifie~:
\begin{eqnarray*}
 \frac12 \Big\{ \frac{\Gamma'}{\Gamma}\Big(\frac{s+\mathfrak{a}}{2}\Big)
+\frac{\Gamma'}{\Gamma}\Big(\frac{1-s+\mathfrak{a}}{2}\Big) \Big\} = \Re
\frac{\Gamma'}{\Gamma}\Big(\frac{s+\mathfrak{a}}{2}\Big)
\end{eqnarray*}
Gr\^ace \`a la condition (B), nous obtenons~:
\begin{multline}
\label{3.5}
\lim_{T\rightarrow +\infty} I_3(T)
=\frac1{2i\pi} \lim_{T\rightarrow+\infty} \int_{1/2-iT}^{1/2+iT} \Re
\frac{\Gamma'}{\Gamma}\Big(\frac{s+\mathfrak{a}}{2}\Big) \Phi(-s) d\,s
\\ + (1-\mathfrak{a}) \Phi(0)
\end{multline}
Et \myref{3.1}, \myref{3.2}, \myref{3.3}, \myref{3.4} et \myref{3.5}
donnent ainsi l'\'egalit\'e annonc\'ee.
\end{proof}
\noindent Pour obtenir \myref{1.1}, il ne reste plus qu'\`a prendre la partie
r\'eelle dans la formule du th\'eor\`eme \ref{FormuleExplicite}
dans le cas $q=1$ et $s=\sigma+it$, o\`u a priori $\sigma>1$, et
$$
 \phi(y) = \begin{cases}
(f(0)-f(y))e^{-ys} & \mathrm{si\ } y\ge0 \\
0 & \mathrm{sinon}
\end{cases}
$$
$\phi$ est alors une fonction de classe $C^2$ sur $\mathbb{R}$ qui v\'erifie
bien la condition (B) et on a pour $\Re z < \Re s$~:
\begin{gather*}
\Phi(-z) = \frac{f(0)}{s-z} - F(s-z) =-\frac{F_2(s-z)}{(s-z)^2},
\\ \Phi(0) =-\frac{F_2(s)}{s^2}, \quad \Phi(-1) = \frac{f(0)}{s-1} - F(s-1)
\end{gather*}
o\`u $F_2$ est la transform\'ee de Laplace de $f"$.
Il vient  alors
\begin{multline}
\label{2.2}
\Re \sum_{n\ge1} \frac{\Lambda(n)}{n^{s}} f(\log n)
= f(0) \ \Re \Bigg( \sum_{n\ge1} \frac{\Lambda(n)}{n^{s}} - \frac1{s-1} +
\sum_{\varrho\in Z(\zeta)}\frac1{s-\varrho} \Bigg) \\
+ \Re F(s-1) -
\sum_{\varrho\in Z(\zeta)} \Re F(s-\varrho)
\\ +
\Re \Bigg( \frac{1}{2i\pi}\int_{1/2-i\infty}^{1/2+i\infty}
\Re\frac{\Gamma'}{\Gamma}\Big(\frac z2 \Big)
\frac{F_2(s-z)}{(s-z)^2}d\,z + \frac{F_2(s)}{s^2} \Bigg)
\end{multline}
La formule d'Hadamard (voir \cite{Dav}) permet de r\'e\'ecrire le terme facteur de
$f(0)$ ~:
$$
\sum_{n\ge1} \frac{\Lambda(n)}{n^s}
= -\frac{\zeta'}{\zeta} (s)
= -B - \frac12\log \pi + \frac1{s-1} + \frac12 \frac{\Gamma'}{\Gamma}\Big(\frac{s}2+1\Big)
- \sum_{\varrho \in Z(\zeta)}\Big(\frac1{\varrho}+\frac1{s-\varrho}\Big)
$$
or $\Re B = -\sum_{\varrho \in Z(\zeta)} \Re \frac1{\varrho}$ donc
$$
 \Re \Big( -\frac{\zeta'}{\zeta}(s) - \frac1{s-1} + \sum_{\varrho\in
Z(\zeta)}\frac1{s-\varrho} \Big)
= -\frac12\log \pi + \frac12\Re
\frac{\Gamma'}{\Gamma}\Big(\frac{s}2+1\Big)
$$
L'identit\'e \myref{2.2} devient~:
\begin{multline*}
\Re \sum_{n\ge1} \frac{\Lambda(n)}{n^{s}} f(\log n)
= f(0) \Bigg( - \frac12 \log(\pi) + \frac12
\Re \frac{\Gamma'}{\Gamma}\Big(\frac{s}2+1\Big) \Bigg)
\\ + \Re F(s-1) - \sum_{\varrho\in Z(\zeta)} \Re F(s-\varrho)
\\ + \Re \Bigg( \frac{1}{2i\pi}\int_{1/2-i\infty}^{1/2+i\infty}
\Re\frac{\Gamma'}{\Gamma}\Big(\frac z2 \Big)
\frac{F_2(s-z)}{(s-z)^2}d\,z + \frac{F_2(s)}{s^2} \Bigg)
\end{multline*}
Nous avons ici une \'egalit\'e entre deux fonctions harmoniques sur le
demi-plan $\Re s > 1$, mais les deux membres d\'efinissent des
fonctions sur $\mathbb{C}$ tout entier
(l'introduction de la partie r\'eelle \^ote les probl\`emes de convergence)
ce qui fait que l'\'egalit\'e reste vraie sur $\mathbb{C}$.
Ceci ach\`eve la d\'emonstration de la proposition \ref{Prop1.1}.
\subsection{\'Etude de $ \tilde F$.}
\label{Section5.2}\noindent
Nous \'etudions dans ce paragraphe le comportement de la fonction
$\tilde F$ qui, rappelons-le, d\'epend des param\`etres $\theta$ et
$\eta$~:
$$
\tilde F (x,y) = \int_0^{d_1(\theta)/\eta} e^{-xt} \cos (yt) f(t) d\,t
 = \int_0^{d_1(\theta)} \exp \Big(\frac{-xt}{\eta}\Big) \cos
 \Big(\frac{yt}{\eta}\Big) h_{\theta} (t) d\,t.
$$
Nous avons fix\'e $\theta=1.848$ et les r\'esultats num\'eriques
donn\'es ici sont calcul\'es pour cette valeur.
\begin{lem}
\label{Lemme2.2} Nous avons
\begin{equation*}
\tilde F(x,y) = \eta g_1(\theta) \frac{x}{x^2+y^2} +H(x,y)
\end{equation*}
o\`u la fonction $H$ v\'erifie
\begin{equation*}
|H(x,y)| \le \frac{M( x /\eta)\eta^2}{x^2+y^2},
\quad\hbox{avec}\quad M(z) =\int_0^{d_1(\theta)}|h_{\theta}^"(u)|e^{-zu}d\,u.
\end{equation*}
De plus
lorsque $0 \le z \le \frac1{d_1(\theta)}$, nous avons le d\'eveloppement suivant~:
$$
521.632 -212.574 z\le
M(z) \le 521.633 - 212.573 z + 68.114 z^2,
$$
Sinon, la majoration par $m/z$ donne
une approximation de $M$ suffisante o\`u
$$
m=\max_{u\in[0,d_1(\theta)]}|h_{\theta}^"(u)| = |h_{\theta}^"(0)|  =
1322.86625 + \mathcal{O^*}(10^{-5}).
$$
\end{lem}
\begin{proof} \ \newline\noindent
Rappelons que~:
$$
F(s) = \frac{f(0)}{s}+\frac{F_2(s)}{s^2}
$$
o\`u $F_2$ est la transform\'ee de laplace de $f^"$.  Donc~:
\begin{eqnarray*}
\tilde F (x,y)&=& \Re \Big[ \frac{f (0)}{x+iy}
+\frac1{(x+iy)^2}\int _0^{d(\theta,\eta)} e^{-(x+iy)t} f^" (t) \,dt\Big]
\\&=&f (0)\frac{x}{x^2+y^2} +
\int _0^{d(\theta,\eta)} \frac{(x^2-y^2)\cos(ty)-2xy\sin(ty)}{(x^2+y^2)^2}
e^{-xt}f^" (t) \,dt
\end{eqnarray*}
Notons $H$ le reste et majorons le :
\begin{equation}
\label{DefH}
H(x,y) = \int _0^{d(\theta,\eta)}
\frac{(x^2-y^2)\cos(ty)-2xy\sin(ty)}{(x^2+y^2)^2}e^{-xt}f^" (t) \,dt.
\end{equation}
L'in\'egalit\'e de Cauchy-Schwartz nous donne
\begin{equation*}
\label{firsttrick}
|(x^2-y^2)\cos(ty)-2xy\sin(ty)|\le
\sqrt{(x^2-y^2)^2+(2xy)^2}=x^2+y^2
\end{equation*}
et par cons\'equent
$$
|H(x,y)| \le  \frac{1}{x^2+y^2}
\int_0^{d(\theta,\eta)}e^{-xt} |f^" (t)| \,dt
$$
Par d\'efinition de $f$, $f^"(t) = \eta^3h_{\theta}^"(\eta t)$, donc un
changement de variable donne que l'int\'egrale de droite est \'egale
\`a $\eta^2\int_0^{d_1(\theta)}|h_{\theta}^"(u)|e^{-\frac{x}{\eta}u}d\,u$.
\par \noindent
De plus, les in\'egalit\'es \'el\'ementaires $1-t\le e^{-t}\le
1-t+\frac{t^2}2$, valables pour $t\ge-1$, nous donnent l'encadrement annonc\'e
pour $M(z)$, ce
qui ach\`eve la d\'emonstration du lemme.
\end{proof}
\noindent Nous aurons besoin au paragraphe \ref{Section2.5} d'une estimation
plus pr\'ecise de $H(x,y)$ lorsque y tend vers l'infini~:
\begin{lem}
\label{Lemme2.4}
Nous avons
\begin{eqnarray*}
|H(x,y)| \le m\eta^3 \frac{|x||x^2-3y^2|}{(x^2+y^2)^3} +
\frac{M_1(x/\eta)\eta^3}{(x^2+y^2)^{3/2}}
\\ \hbox{avec} \qquad M_1(z) =\int_0^{d_1(\theta)}|h_{\theta}^{(3)}(u)|e^{-zu}d\,u.
\end{eqnarray*}
De plus
lorsque $0 \le z \le \frac1{d_1(\theta)}$, nous avons le d\'eveloppement suivant~:
$$
2526.445 -1087.743 z \le
M_1(z) \le 2526.446 - 1087.742 z + 348.808 z^2,
$$
Sinon, la majoration par $m_1/z$ donne
une approximation de $M_1$ suffisante avec
$$
m_1=\max_{u\in[0,d_1(\theta)]}|h_{\theta}^{(3)}(u)| = 4135.12706 + \mathcal{O}^*(10^{-5}).
$$
\end{lem}
\begin{proof} \ \newline\noindent
En int\'egrant par parties
$$
H(x,y) = \Re \Big( \frac1{(x+iy)^2} \int_0^{d(\theta,\eta)}e^{-(x+iy)t}f^" (t)d\,t \Big)
$$
nous obtenons
\begin{eqnarray*}
H(x,y)&=&\Re \Big( \frac{f^" (0)}{(x+iy)^3} + \frac1{(x+iy)^3}\int_0^{d(\theta,\eta)}e^{-(x+iy)t}f^{(3)}(t)d\,t \Big)
\\ &=&f^" (0) \frac{x(x^2-3y^2)}{(x^2+y^2)^3}
\\ \quad&&
+\int_0^{d(\theta,\eta)} \frac{x(x^2-3y^2)\cos(yt)-y(y^2-3x^2)\sin(yt)}{(x^2+y^2)^3} e^{-xt} f^{(3)}(t)d\,t
\end{eqnarray*}
En achevant la preuve comme celle du lemme \ref{Lemme2.2}, nous montrons ainsi que~:
$$
H(x,y) \le h_{\theta}^" (0)\eta^3\frac{|x||x^2-3y^2|}{(x^2+y^2)^3} +
\frac{\eta^3}{(x^2+y^2)^{3/2}} \int_0^{d_1(\theta)}|h_{\theta}^{(3)(u)}|
e^{-\frac{x}{\eta}u} d\,u
$$
\end{proof}
\noindent Dor\'enavant, nous op\'erons \`a $y\ge0$ fix\'e.  D'apr\`es le
r\'esultat ci-dessus, nous pouvons esp\'erer que , pour $x \in
[0,+\infty[$, la fonction $x \mapsto \tilde F( x, y)$ se
comporte comme la fonction
$x \mapsto \frac{x}{x^2+y^2}$, c'est \`a dire qu'elle croisse jusqu' \`a une valeur
proche de y puis d\'ecroisse ensuite.
\noindent
Lorsque $y=0$, il est imm\'ediat de voir que la fonction est
d\'ecroissante et tend vers 0 en l'infini.
\noindent
Lorsque y est strictement positif, on a le lemme suivant:
\begin{lem}
\label{Lemme2.3}
Pour tout $y>0$ l'application $\tilde F(\cdot, y)$
d\'ecro\^\i t
sur $[x_2(\eta,y),+\infty[$ avec
$x_2(\eta,y) =\varepsilon_3(\eta)+\sqrt{y^2+\varepsilon_3^2(\eta)}$,
$\varepsilon_3(\eta)=7.857\eta$.
\noindent
Par ailleurs, pour tout $y>0$ l'application $\tilde F(\cdot, y)$
cro\^\i t
sur $[0,x_1(\eta,y)]$ avec
$x_1(\eta,y)=\tfrac12y-\varepsilon_1(\eta)+\sqrt{(\tfrac12y-\varepsilon_1(\eta))^2-
\varepsilon_2(\eta)}$ pourvu que la quantit\'e sous la racine soit
positive ou nulle, o\`u nous avons pos\'e
$\varepsilon_1(\eta)=4.99\eta$ et $\varepsilon_2(\eta)=5.735\eta$.
\end{lem}
\begin{proof} \ \newline\noindent
Reprenons $H$ de la d\'emonstration pr\'ec\'edente
donn\'e par \myref{DefH}.
Il vient
\begin{multline*}
\frac{\partial}{\partial x}H( x, y) =
\int _0^{d(\theta,\eta)} \Bigg( -t\frac{(x^2-y^2)\cos(t y)-2xy\sin(t y)}{(x^2+y^2)^2}
\\ \qquad -2\frac{x(x^2-3y^2)\cos(t y)-y(3x^2-y^2)\sin(t y)}{(x^2+y^2)^3}
\Bigg) e^{- xt} f" (t) \,dt.
\end{multline*}
Nous utilisons encore l'in\'egalit\'e de Cauchy-Schwartz pour montrer que~:
$$
|2x(x^2-3y^2)\cos(ty)-2y(3x^2-y^2)\sin(ty)|\le
2(x^2+y^2)^{3/2}
$$
ainsi que le lemme \ref{Lemme2.2} pour obtenir
\begin{eqnarray*}
\bigg|\frac{\partial}{\partial x}H( x, y)\bigg| &\le&
\int _0^{d(\theta,\eta)} |f"(t)| \Bigg( \frac{t} {x^2+y^2}
+\frac{1}{(x^2+y^2)^{3/2}} \Bigg) e^{-xt} \,dt \\
&& \quad =\eta \frac{M_2(x/\eta)}{x^2+y^2}
+ 2\eta^2 \frac{M(x/\eta)}{(x^2+y^2)^{3/2}}
\end{eqnarray*}
$$\hbox{o\`u}\quad M_2(z) = \int_0^{d_1(\theta)}|h_{\theta}^"(u)|ue^{-zu}d\,u
\le \|u h''_\theta\|_\infty/z.
$$
Supposons $x\ge y$. Puisque
\begin{equation}
\label{partialF}
\frac{\partial}{\partial x}\tilde F( x, y) = g_1(\theta)\eta
\frac{y^2-x^2}{(x^2+y^2)^2} +\frac{\partial}{\partial x}H( x, y)
\end{equation}
nous pouvons garantir que  $\frac{\partial}{\partial x}F( x, y)\le0$ d\`es
que
$$
g_1(\theta)(y^2-x^2)
+M_2(x/\eta)(x^2+y^2)+2\eta M(x/\eta)\sqrt{x^2+y^2}
\le 0
$$
ce qui est impliqu\'e par
$$
g_1(\theta)(y^2-x^2)
+2\eta(\|uh''_\theta\|_\infty +\sqrt{2}\|h''_\theta\|_1) x
\le 0
$$
avec $\|uh''_\theta\|_\infty\le423.867$ et $\|h''_\theta\|_1\le521.633$.
$$
x\ge x_2(\eta,y)=
\varepsilon_3(\eta)+\sqrt{y^2+\varepsilon_3^2(\eta)}
$$
$$\hbox{avec}\quad
\varepsilon_3(\eta)=7.857\eta\ge
\frac{\|uh''_\theta\|_\infty
  +\sqrt{2}\|h''_\theta\|_1}{g_1(\theta)}\eta.
$$
\`A partir de \myref{partialF}, nous pouvons aussi garantir que
$\frac{\partial}{\partial x}F( x, y)\ge0$ d\`es que
$x\le y$
$$\hbox{et}\quad
g_1(\theta)(y^2-x^2)
-M_2(x/\eta)(x^2+y^2)-2\eta M(x/\eta)\sqrt{x^2+y^2}
\ge 0
$$
ce qui est impliqu\'e par
$$
g_1(\theta)(y^2-x^2)
-\eta\|uh''_\theta\|_\infty\frac{2y^2}{x} -2\sqrt{2}\eta\|h''_\theta\|_1y
\ge 0.
$$
Comme nous n'aurons pas besoin d'un r\'esultat tr\`es performant, nous
nous contentons de noter que $y^2-xy\le y^2-x^2$, ce qui nous laisse
avec
$$
g_1(\theta)(y-x)x
-2\eta\|uh''_\theta\|_\infty y -2\sqrt2\eta\|h''_\theta\|_1x
\ge 0
$$
et il nous suffit maintenant d'avoir
$$
x^2-(y- 9.980\eta)x+5.735\eta y\le0.
$$
\end{proof}
\subsection{\'Etude de $\Re \frac{\Gamma'}{\Gamma}$}
\label{Section5.3}\noindent
Dans la suite nous allons avoir besoin d'une estimation du terme $\Re
\frac{\Gamma'}{\Gamma}$ pour \'etudier les termes $\Delta_1$ et $\Delta_2$. C'est l'objet des deux lemmes suivants.
\noindent
\begin{lem}
\label{Gamma1}
Soient $\delta \in [0,1]$, $\kappa \in [0,x/(x+\delta)]$ et $0<x_0\le
x \le x_1 < y_0$, alors~:
\begin{multline*}
\Re\frac{\Gamma'}{\Gamma}\Big(\frac{x}2+i\frac{y}2\Big) - \kappa\ \Re\frac{\Gamma'}{\Gamma}\Big(\frac{x+\delta}2+i\frac{y}2\Big)\\ \le
\begin{cases}
r_1(x_0,x_1,y_0) \ \mathrm{si}\ 0<|y|< y_0\\
(1-\kappa)\log \frac{y}2 + \min \Big(r_2(x_0,x_1,y_0),r_3(x_0,x_1,y_0)\Big)\ \mathrm{si}\ |y|\ge y_0 \\
\end{cases}
\end{multline*}
où
$r_1$, $r_2$ et $r_3$ sont respectivement d\'efinis en \myref{r1}, \myref{r2}, \myref{r3}.
\end{lem}
\begin{proof} \ \newline\noindent
Pour la suite, notons $\psi_{\kappa,\delta}(x,y)$ la diff\'erence $\Re\frac{\Gamma'}{\Gamma}\Big(\frac{x}2+i\frac{y}2\Big) - \kappa\ \Re\frac{\Gamma'}{\Gamma}\Big(\frac{x+\delta}2+i\frac{y}2\Big)$.
Nous allons approcher le terme
$\Re\frac{\Gamma'}{\Gamma}\Big(\frac{x}2+i\frac{y}2\Big)$ de deux fa\c cons
diff\'erentes qui sont plus ou moins efficaces selon la taille de $|y|$.
\par\noindent
Utilisons l'identit\'e donn\'ee par K.Mc.Curley (Cf. \cite{KSMC})~:
\begin{equation}
\label{gamma1}
\Re\frac{\Gamma'}{\Gamma}\Big(\frac{x}2+i\frac{y}2\Big)
= \frac12 \log\Big(\frac{x^2}4+\frac{y^2}4\Big) - \frac{x}{x^2+y^2}
+ \Re\int_{0}^{+\infty} \frac{(u-[u]-1/2)}{\big(u+\frac{x+iy}2\big)^2} d\,u
\end{equation}
avec le terme int\'egral qui satisfait~:
$$
\Re\int_{0}^{+\infty}
\Bigg|\frac{(u-[u]-1/2)}{\Big(u+\frac{x+iy}2\Big)^2}\Bigg|d\,u \le \frac1{y} \arctan \frac{y}{x}
$$
\myref{gamma1} permet ainsi de majorer $\psi_{\kappa,\delta}(x,y)$ par $R_1(x,y)$ d\'efini par~:
\begin{multline}
 \frac12 \log\Big(\frac{x^2}4+\frac{y^2}4\Big)
 - \frac{\kappa}2 \log\Big(\frac{(x+\delta)^2}4+\frac{y^2}4\Big)
 \\- \Big(\frac{x}{x^2+y^2}-\kappa\ \frac{x+\delta}{(x+\delta)^2+y^2}\Big)
 + \frac1{y}\Big(\arctan \frac{y}{x}+\kappa\arctan \frac{y}{x+\delta}\Big)
\end{multline}
Nous allons maintenant distinguer les cas o\`u $y$ est born\'e ($0<|y|<y_0$) et o\`u $y$ est ``grand'' ($|y|\ge y_0$).
Dans le premier cas, nous majorons $\frac1{y} \arctan \frac{y}x$ par $\frac1{x}$ et nous obtenons
\begin{multline}
\label{r1}
R_1(x,y) \le r_1(x_0,x_1,y_0) = \frac{1-\kappa}2 \log \Big(\frac{(x_1+\delta)^2}4+\frac{y_0^2}4\Big)
\\- \frac{x_0}{x_1^2+y_0^2} + \frac1{x_0} + \frac{2\kappa}{x_0+\delta}\quad \hbox{si}\ 0<|y|<y_0
\end{multline}
Dans le second cas, $\Re\frac{\Gamma'}{\Gamma}\Big(\frac{x}2+i\frac{y}2\Big)$ peut être approch\'e par $\log |y|$. Nous r\'e\'ecrivons $R_1(x,y)$~:
$$
R_1(x,y) = (1-\kappa)\log\frac{|y|}2 + R_2(x,y)
$$
\begin{multline}
\label{R2}
\hbox{avec}\quad R_2(x,y) = \frac12\log\Big(\frac{x^2}{y^2}+1\Big) - \frac{\kappa}2 \log \Big(\frac{(x+\delta)^2}{y^2}+1\Big)
\\- \Big(\frac{x}{x^2+y^2}
  -\kappa\ \frac{x+\delta}{(x+\delta)^2+y^2}\Big)
+\frac1{y}\Big(\arctan \frac{y}{x}+\kappa\ \arctan \frac{y}{x+\delta}\Big)
\end{multline}
Comme l'application $y \mapsto \frac1{y}\arctan \frac{y}{x}$
est positive et d\'ecroissante ainsi que $y \mapsto \frac{x}{x^2+y^2}-\kappa\
\frac{x+\delta}{(x+\delta)^2+y^2}$
puisque $\kappa\le \frac{x}{x+\delta}$, nous obtenons pour le terme d'erreur~:
\begin{multline}
\label{r2}
R_2(x,y) \le  r_2(x_0,x_1,y_0) = \frac{1-\kappa}2 \log
\Big(\frac{(x_1+\delta)^2}{y_0^2}+1\Big) + \frac1{y_0}\Big(\arctan \frac{y_0}{x_1}
\\+\kappa\ \arctan \frac{y_0}{x_1+\delta}\Big)\quad\hbox{si}\ |y|\ge y_0
\end{multline}
L'\'egalit\'e suivante donne une autre majoration pour le terme d'erreur~:
\begin{equation}
\label{Gamma3}
\Re\frac{\Gamma'}{\Gamma}(x+iy) = \log|y| - \frac{x}{2(x^2+y^2)} + R(x,y) \qquad  (y^2>x^2,x>0)
\end{equation}
$$\hbox{avec}\quad
|R(x,y)| \le \frac{1}{12x|y|}+\frac{x^2}{2y^2}
$$
Et donc ~:
$$
\psi_{\kappa,\delta}(x,y)
\le (1-\kappa) \log \frac{|y|}{2} + R_3(x,y)
 \le (1-\kappa) \log \frac{|y|}{2} + r_3(x_0,x_1,y_0)
$$
\begin{multline}
\label{R3}
\hbox{avec}\quad
R_3(x,y) = - \Big(\frac{x}{x^2+y^2}-\kappa\frac{x+\delta}{(x+\delta)^2+y^2}\Big) + \frac1{3y}\Big(\frac1{x}+\frac{\kappa}{x+\delta}\Big)
\\ + \frac1{2y^2}\Big(x^2+\kappa(x+\delta)^2\Big)
\end{multline}
\begin{equation}
\label{r3}
r_3(x_0,x_1,y_0) = \frac1{3y_0}\Big(\frac1{x_0}+\frac{\kappa}{x_0+\delta}\Big) + \frac1{2y_0^2}\Big(x_1^2+\kappa(x_1+\delta)^2\Big)
\end{equation}
\end{proof}
\begin{lem}
\label{Gamma2}
Nous avons~:
$$
\Big|\Re\frac{\Gamma'}{\Gamma}\Big(\frac14+i\frac{T}{2}\Big)\Big|\le
U_0(T) =
\begin{cases}
\frac12\log\frac{16}{1+4T^2} + \frac2{1+4T^2} - \frac{\pi}{2} \ \mathrm{si}\ |T|<1/2 \\
\big|\log \frac{|T|}2 - \frac2{1+4T^2} \big| + \frac2{3|T|} + \frac1{8T^2} \ \mathrm{si}\ |T| \ge 1/2
\end{cases}
$$
\end{lem}
\begin{proof} \ \newline\noindent
Pour la suite, nous d\'esignerons par  $R_i(T)$ l'application $T \mapsto R_i(1/2,T,0,0)$.
\begin{enumerate}
\item
Si $|T|<1/2$, nous avons grâce au lemme \ref{Gamma1}~:
\begin{equation}
\label{in1}
\Bigg| \Re \frac{\Gamma'}{\Gamma}\Big(\frac{1}4+i\frac{T}2\Big) \Bigg|
\le |R_1(T)|
\end{equation}
avec
$$
R_1(T) = \frac12 \log \Big(\frac1{16}+\frac{T^2}4\Big) -
\frac{2}{1+4T^2}
+ \frac1{T} \arctan(2T)
$$
$R_1$ est n\'egative sur $[0,1/2]$, donc \myref{in1} devient
\begin{equation}
\Big|\Re \frac{\Gamma'}{\Gamma}\Big(\frac14+i\frac{T}2\Big)\Big|
\le -R_1(T)
 \le \frac12 \log \Big(\frac{16}{1+4T^2}\Big) + \frac2{1+4T^2} - \frac{\pi}2
\end{equation}
\item
Si $|T|\ge1/2$, l'\'egalit\'e \myref{Gamma3} nous donne~:
\begin{equation*}
\Re\frac{\Gamma'}{\Gamma}\Big(\frac14+i\frac{T}2\Big)
 = \log \frac{|T|}2 - \frac2{1+4T^2} + R\Big(\frac14,\frac{|T|}2\Big)
\end{equation*}
et donc
\begin{equation}
\label{in2}
\Big|\Re\frac{\Gamma'}{\Gamma}\Big(\frac{1}4+i\frac{T}2\Big)\Big|
\le \Big| \log \frac{|T|}2 - \frac2{1+4T^2}\Big| + \frac2{3|T|} + \frac{1}{8T^2}
\end{equation}
\end{enumerate}
\end{proof}
\section{Preuves}
\label{Section6}\noindent
Dans ce paragraphe, nous commen\c cons par \'etudier la somme sur les
z\'eros de Z\^eta. Le r\'esultat fondamental de la proposition
\ref{Prop2.6} permet de r\'egler le cas des z\'eros de partie r\'eelle
parcourant $[1-\sigma,\sigma]$ par un argument de positivit\'e.\\
\noindent
Nous introduisons \`a cette occasion les conditions num\'eriques
suivantes sur les param\`etres $\kappa$ et $\delta$ ~:
\begin{gather*}
0 \le \kappa \le 0.4389 \quad \hbox{ et } \quad \delta \ge 0.62063,
\end{gather*}
ce qui nous permet d'obtenir les approximations num\'eriques annonc\'ees
aux propositions \ref{Prop1.2}, \ref{Prop1.3}, \ref{Prop1.4} et
\ref{Prop1.5}. \\ \noindent
Nous rappelons \`a cette occasion les valeurs que nous avons fix\'ees
pour les param\`etres de d\'epart $T_0$, $R$, $\theta$ et $r$~:
\begin{gather*}
T_0 = 3 330 657 430 .697, \quad R=9.645908801,\\
 \quad \theta = 1.848, \quad r=5.97484,
\end{gather*}
ainsi que celles que nous obtenons pour les variables interm\'ediaires
dont nous avons besoin ici ~:
\begin{gather*}
\quad g_1(\theta)=147.84112 , \quad \sigma_0 =0.99555 ,\quad \eta_0=0.00913 , \\
\quad m=1322.86625 ,\quad m_1=4135.12706,\\
M(0)= 521.632466 ,\quad M(-1)= 822.67426.
\end{gather*}
\subsection{\'Etude de la somme sur les z\'eros }
\label{Section2.5} \noindent
Nous cherchons \`a localiser le z\'ero $\varrho_0 = \beta_0 + i \gamma_0$.
Pour cela, et c'est l'objet de ce premier paragraphe, nous l'isolons dans la somme
\begin{equation*}
\sum_{k=0}^4 a_k \sum_{\varrho \in Z(\zeta)}
D(\sigma+ik\gamma_0-\varrho).
\end{equation*}
\subsubsection{Cas $k=1$ et $\varrho \in
  \{\varrho_0,1-\overline\varrho_0\}$~:}\noindent
Le terme $D(\sigma-\beta_0) + D(\sigma-1+\beta_0)$ donn\'e par
\begin{eqnarray*}
D(\sigma-\beta_0) &=& \tilde F(\sigma-\beta_0,0) -\ \kappa\ \tilde
F(\sigma-\beta_0 +\delta,0)
\\ \hbox{et}\quad
D(\sigma-1+\beta_0) &=& \tilde F(\sigma-1+\beta_0,0) -\ \kappa\ \tilde
F(\sigma-1+\beta_0 +\delta,0)
\end{eqnarray*}
est en fait proche de $\tilde F (\sigma-\beta_0,0)$ \` a un
$\mathcal{O}(\eta)$ pr\`es.\\ \\ \noindent
En effet, d'apr\`es la  d\'ecroissance de l'application $x \mapsto \tilde F(x,0)$ et le lemme \ref{Lemme2.2}, nous avons~:
\begin{eqnarray*}
\tilde F(\sigma-1+\beta_0 +\delta,0) \le \tilde F(1-\eta_0 +\delta,0)
&\le& \frac{g_1(\theta)}{1-\eta_0+\delta}\eta + \frac{m}{(1-\eta_0+\delta)^3}\eta^3
\\ \tilde F(\sigma-\beta_0+\delta,0) &\le& \tilde F(\delta,0) \ge \frac{g_1(\theta)}{\delta} \eta +
\frac{m}{\delta^3} \eta^3
\\ \tilde F(\sigma-1+\beta_0,0) &\ge& \tilde F(1,0) \ge g_1(\theta) \eta - m \eta^3
\end{eqnarray*}
Et donc~:
\begin{multline}
\label{3.10}
 D(\sigma-\beta_0) + D(\sigma-1+\beta_0) \ge \tilde F(\sigma-\beta_0,0)
\\+\Big(g_1(\theta)\eta-m\eta^3\Big)
-\ \kappa \Bigg(\Big(\frac1{\delta}+\frac1{1-\eta_0+\delta}\Big)g_1(\theta)\eta +
\Big(\frac1{\delta^3}+\frac1{(1-\eta_0+\delta)^3}\Big) m \eta^3 \Bigg)
\end{multline}
Nous \'etudions maintenant le reste de la somme et nous allons montrer,
gr\^ace aux propositions \ref{Prop2.6} et \ref{Prop2.7}, qu'en fait il est d'ordre $\mathcal{O}(\eta^2)$.
\subsubsection{Cas $(k=0,2,3,4)$ ou $(k=1$ et $\rho \not\in \{\rho_0,1-\overline{\rho_0}\})$.}
\noindent Les z\'eros de $\zeta$ \'etant sym\'etriques par rapport
\`a l'axe r\'eel et \`a l'axe $\Re s = 1/2$, nous avons~:
\begin{eqnarray*}
&&\sum_{\varrho \in Z(\zeta)} D (\sigma+ik\gamma_0- \varrho)
= \frac12 \sum_{\varrho\in Z(\zeta)} \Big[
D(\sigma+ik\gamma_0-\varrho)+D (\sigma-1+ik\gamma_0 + \bar
\varrho)\Big]
\\ && =  \frac12 \sum_{{\begin{substack}{\varrho \in Z(\zeta) \\
                       \beta > 1/2}
      \end{substack}}}
\Big[ D(\sigma -\beta+i(k\gamma_0-\gamma))+D
(\sigma-1+\beta+i(k\gamma_0-\gamma)) \Big]
\\ &&\qquad + \sum_{{\begin{substack}{\varrho \in Z(\zeta) \\
                       \beta = 1/2}
      \end{substack}}}
  D(\sigma -\tfrac12+i(k\gamma_0-\gamma))
\end{eqnarray*}
L'argument de positivit\'e de la proposition suivante nous permet
d'\'eliminer une grande partie des z\'eros. Nous g\'en\'eralisons \`a
la transform\'ee de Laplace $F$ le r\'esultat de Stechkin (voir
\cite{Stechkin1}) ~:
\begin{lem}[Stechkin - 1970]\ \\ \noindent
Pour $\beta \in [ \tfrac12 , 1 ] $, $y>0$, $\sigma >1$ et
$\tau = \frac{1+\sqrt{1+4\sigma^2}}{2}$, nous avons
\begin{equation*}
\Re \Big(\frac1{\sigma-\beta+iy}\, -\, \frac1{\sqrt5}
\frac1{\big(\tau-\beta+iy\big)} \Big) \ +\  \Re \Big(\frac1{\sigma-1+\beta+iy} \,-\,
\frac1{\sqrt5} \frac1{\big(\tau-1+\beta+iy\big)} \Big) \ge 0.
\end{equation*}
\end{lem}
\begin{prop}\ \\
\label{Prop2.6}
Pour $\beta\in[\tfrac12,\sigma]$ et $y>0$, nous avons
\begin{equation*}
D (\sigma -\beta+iy)+D (\sigma-1+\beta+iy) \ge 0
\end{equation*}
d\`es que\  $0\le\kappa\le 0.4389$\  et\  $\delta \ge 0.62063$.
\end{prop}
\begin{proof} \ \newline\noindent
Nous cherchons le plus grand $\kappa$ tel que
$$
D(\sigma-\beta+iy) + D(\sigma-1+\beta+iy) \ge 0.
$$
C'est \`a dire que nous cherchons \`a minorer la fonction suivante~:
$$
Q(\beta+iy) = \frac{\tilde F(\sigma-\beta,y)+\tilde
  F(\sigma-1+\beta,y)}{\tilde F(\sigma+\delta-\beta,y)+\tilde
  F(\sigma+\delta-1+\beta,y)}
$$
Nous allons tout d'abord montrer que le num\'erateur de $Q$ ne
s'annule jamais sur $\mathbb{C}$, ce qui signifiera que $Q$, en tant
que quotient de parties r\'eelles de fonctions enti\`eres, est une
fonction harmonique sur $\mathbb{C}$~:
la positivit\'e de $\tilde F$ implique que $\tilde
F(\sigma+\delta-\beta,y)+\tilde F(\sigma-1+\delta+\beta,y)$ s'annule
si et seulement si $\tilde F(\sigma+\delta-\beta,y)$ et $\tilde
F(\sigma-1+\delta+\beta,y)$ s'annulent.
En utilisant les majorations du lemme \ref{Lemme2.2}, nous obtenons
les deux in\'egalit\'es~:
\begin{eqnarray*}
&&\frac{g_1(\theta)(\sigma-\beta+\delta)}{(\sigma-\beta+\delta)^2+y^2}-\frac{m\eta_0^2}{(\sigma-\beta+\delta)((\sigma-\beta+\delta)^2+y^2)}
\le 0
\\
&&\frac{g_1(\theta)(\sigma-1+\beta+\delta)}{(\sigma-1+\beta+\delta)^2+y^2}-\frac{m\eta_0^2}{(\sigma-1+\beta+\delta)((\sigma-1+\beta+\delta)^2+y^2)}
\le 0
\end{eqnarray*}
ce qui \'equivaut \`a ce que $|\sigma-\beta+\delta|$ et
$|\sigma-1+\beta+\delta|$ soient tous deux major\'es par
$\sqrt{\frac{m}{g_1(\theta)}}\eta_0$ et donc que $\delta$ soit
major\'e par le terme n\'egatif
$\sqrt{\frac{m}{g_1(\theta)}}\eta_0+1/2-\sigma$, ce qui est absurde.

Nous pouvons maintenant appliquer le principe du maximum \`a $Q$~:
nous constatons d'une part qu'il suffit de chercher son minimum sur un
contour pour l'obtenir \`a l'int\'erieur du dit contour, d'autre part
que l'application $y_0 \mapsto \min_{1-\sigma\le\beta\le\sigma,|y|\le
  y_0} Q(\beta+iy)$ est d\'ecroissante sur $[0,+\infty[$. Nous pouvons
donc supposer $y_0$ assez grand, au moins  sup\'erieur \`a $\sigma
-\beta$ et $\sigma-1+\delta$ (nous prendrons $y_0 \ge 10$).
Par ailleurs, en prenant pour domaine $|y| \le y_0$ et $1-\sigma \le
\beta \le \sigma$, comme $Q(1-z)=Q(z)$ et $Q(\overline{z})=Q(z)$, il
nous suffit de nous restreindre aux deux côt\'es~:
\ $(y=y_0,\ 1/2 \le \beta \le \sigma) \text{\ et\ } (0 \le y \le y_0,\ \beta=\sigma)$.
\noindent
Dans le premier cas, nous minorons $Q(\beta+iy_0)$ par le terme $Q_1(\sigma,\beta,y_0)$~:
$$
\frac{ \eta g_1(\theta) \Big(
\frac{(\sigma-\beta)}{(\sigma-\beta)^2+y_0^2}
  + \frac{(\sigma-1+\beta)}{(\sigma-1+\beta)^2+y_0^2}\Big)
- H(\sigma-\beta,y_0)- H(\sigma-1+\beta,y_0) } { \eta g_1(\theta)
\Big(
  \frac{(\sigma-\beta+\delta)}{(\sigma-\beta+\delta)^2+y_0^2} +
  \frac{(\sigma-1+\beta+\delta)}{(\sigma-1+\beta+\delta)^2+y_0^2}\Big)
+ H(\sigma-\beta+\delta,y_0)+ H(\sigma-1+\beta+\delta,y_0) }
$$
et nous minorons chaque valeur de $H$ grâce au lemme \ref{Lemme2.4}~:
\begin{eqnarray*}
H(\sigma-\beta,y_0) &\le&
\frac{3my_0^3\eta^4}{((\sigma-\beta)^2+y_0^2)^3} +
\frac{M_1(0)\eta^3}{((\sigma-\beta)^2+y_0^2)^{3/2}}
\\ &\le& (3m\eta_0+M_1(0))\frac{\eta_0^2}{y_0^3} \eta
\\ H(\sigma-1+\beta,y_0) &\le&
\frac{3my_0^3\eta^3(\sigma-1+\beta)}{((\sigma-1+\beta)^2+y_0^2)^3} +
\frac{m_1\eta^4/(\sigma-1+\beta)}{((\sigma-1+\beta)^2+y_0^2)^{3/2}}
\\ &\le& \Big(3m
+\frac{m_1\eta_0}{\sigma_0-1/2}\Big)\frac{\eta_0^2}{y_0^3}\eta
\\ H(\sigma-\beta+\delta,y_0) &\le&
\frac{3my_0^3\eta^3(\sigma-\beta+\delta)}{((\sigma-\beta+\delta)^2+y_0^2)^3} +
\frac{m_1\eta^4/(\sigma-\beta+\delta)}{((\sigma-\beta+\delta)^2+y_0^2)^{3/2}}
\\ &\le&
\Big(3m+\frac{m_1\eta_0}{\delta}\Big)\frac{\eta_0^2}{y_0^3}\eta
\\ H(\sigma-1+\beta+\delta,y_0) &\le&
\frac{3my_0^3\eta^3(\sigma-1+\beta+\delta)}{((\sigma-1+\beta+\delta)^2+y_0^2)^3} +
\frac{m_1\eta^4/(\sigma-1+\beta+\delta)}{((\sigma-1+\beta+\delta)^2+y_0^2)^{3/2}}
\\ &\le&
\Big(3m+\frac{m_1\eta_0}{\delta}\Big)\frac{\eta_0^2}{y_0^3}\eta
\end{eqnarray*}
nous avons ainsi~:
$$
Q_1(\sigma,\beta,y_0) \ge \frac{
  g_1(\theta)(2\sigma-1)\frac{y_0^2}{y_0^2+1}-\Big[
  (3m+3m\eta_0+M_1(0))\eta_0^2+\frac{m_1}{1/2-\eta_0}\eta_0^3\Big]\frac1{y_0}}
{
g_1(\theta)(2\sigma+2\delta-1)+\Big[6m\eta_0^2+\frac{2m_1}{\delta}\eta_0^3\Big]\frac1{y_0}}
$$
Nous voyons facilement que le terme de droite est une fonction
croissante en la variable $\sigma$, donc on peut la minorer par sa
valeur en $\sigma_0$, valeur que nous noterons
$\kappa_1(y_0,\delta)$.
\noindent
Regardons maintenant $Q$ sur l'autre côt\'e $0 \le y \le y_0\ ,\ \beta=\sigma$.
Nous allons utiliser le lemme \ref{Lemme2.2} pour $\tilde
F(2\sigma-1,y)$, $\tilde F(\delta,y)$ et $\tilde
F(2\sigma-1+\delta,y)$.
Pour $\tilde F(0,y)$, le lemme ne suffit plus lorsque $y$ est proche
de 0. A la place, nous utilisons la positivit\'e de $\tilde F$ et nous
minorons ainsi $Q(\sigma+iy)$ par~:
$$
 \frac1{(2\sigma-1)^2+y^2} \
\frac{g_1(\theta)(2\sigma-1)-m\eta_0^2/(2\sigma-1)}{ \frac{\delta
    g_1(\theta)+m\eta_0^2/\delta}{\delta^2+y^2}
  +\frac{(2\sigma-1+\delta)g_1(\theta)+m\eta_0^2/(2\sigma-1+\delta)}{(2\sigma-1+\delta)^2+y^2} }
$$
puis par
$$
Q_2(y) = \frac{1}{1+y^2} \
\frac{g_1(\theta)(1-2\eta_0)-m\eta_0^2/(1-2\eta_0)}{ \frac{\delta
    g_1(\theta)+m\eta_0^2/\delta}{\delta^2+y^2}
  +\frac{(\delta+1)g_1(\theta)+m\eta_0^2/(\delta+1-2\eta_0)}{(\delta+1-2\eta_0)^2+y^2} }
$$
Nous \'etudions le sens de variation de $Q_2$ et pour cela regardons
le signe du d\'enominateur de sa d\'eriv\'ee. A un facteur positif
pr\`es, nous trouvons le trinôme du second degr\'e suivant~:\ $Q_3(y)
= d_2y^2+d_1(\theta)y+d_0$, o\`u les $d_i$ sont des fonctions polynômes de
$\delta$.
\noindent
$d_2$ est positif pour $\delta \le 0.07$, n\'egatif sinon et le
discriminant de $Q_3$ est toujours positif pour $\delta\ge0.03$.
Supposons $\delta>0.07$.
\noindent
Alors $Q_3(y)\ge0$ est successivement positive puis n\'egative sur
$[0,+\infty[$ et donc $Q_2$ est d'abord croissante puis d\'ecroissante
et son minimum est \`a d\'eterminer entre $Q_2(0)$ et $Q_2(y_0)$.
\noindent
En fait nous pouvons même regarder sans trop de perte la limite de
$Q_2$ en l'infini au lieu de $Q_2(y_0)$. Nous noterons respectivement
ces valeurs $\kappa_2(\delta)$ et $\kappa_3(\delta)$~:
\begin{equation}
\label{k2}
 \kappa_2(\delta) =
\frac{g_1(\theta)(1-2\eta_0)-m\eta_0^2/(1-2\eta_0)}
{(1+2\delta)g_1(\theta)+\Big(\frac1{\delta}+\frac1{1+\delta-2\eta_0}\Big)m\eta_0^2}
=  \frac1{1+2\delta} + \mathcal{O}(\eta_0)
\end{equation}
\begin{equation}
\label{k3}
 \kappa_3(\delta) =
\frac{g_1(\theta)(1-2\eta_0)-m\eta_0^2/(1-2\eta_0)}
{\Big(\frac1{\delta}+\frac{1+\delta}{(1+\delta-2\eta_0)^2}\Big)g_1(\theta)
+\Big(\frac1{\delta^3}+\frac1{(1+\delta-2\eta_0)^3}\Big)m\eta_0^2}
 = \frac1{\frac1{\delta}+\frac1{1+\delta}} + \mathcal{O}(\eta_0)
\end{equation}
\noindent
En remarquant que $\kappa_1(.,\delta)$ est une fonction d\'ecroissante
et que $\kappa_2(\delta) \le \kappa_1(10,\delta)$, nous voyons qu'il
ne reste plus pour conclure qu'\`a choisir la valeur optimale de
$\min(\kappa_2(\delta),\kappa_3(\delta))$ lorsque $\delta \in [0,1]$.
A un $\mathcal{O}(\eta_0)$ pr\`es, nous prenons donc $\delta$ tel que ~:
$$
\frac1{1+2\delta} = \frac1{\frac1{\delta}+\frac1{1+\delta}} \quad \text{et donc}\quad \kappa=\frac1{1+2\delta}
$$
C'est \`a dire qu'\`a un $\mathcal{O}(\eta_0)$ pr\`es,
$\delta=\frac{\sqrt5-1}2=0.61803+\mathcal{O^*}(10^{-5})$ et
$\kappa=\frac1{\sqrt5}=0.44721+\mathcal{O^*}(10^{-5})$.
Les calculs exacts donnent~:
$$
\delta=0.62063+\mathcal{O^*}(10^{-5}) \quad \text{et}\quad  \kappa =
0.4389+\mathcal{O^*}(10^{-5}).
$$
\end{proof}
\noindent Dans la suite de cette section, nous posons $y_k= |k T_0-\gamma|$.
\noindent
Il reste \`a \'etudier le cas o\`u $\beta \in [\sigma,1]$, ou plut\^ot
celui o\`u $y_k \ge t_0$, puisque $\sigma \ge
1-\frac1{R\log(kT_0+t_0)}$ et $\beta \le 1-\frac1{R\log\gamma}$.
\noindent
Nous utiliserons le lemme pr\'eliminaire suivant pour montrer la
proposition \ref{Prop2.7} ci-apr\`es. \\
\begin{lem}
\label{Lemme4}
Pour tout $t_0 \ge 1$,
\begin{gather*}\sum_{{\begin{substack}{\varrho \in Z(\zeta) \\
                       |\gamma-t| \ge t_0}
      \end{substack}}} \frac1{(\gamma-t)^2}
\le \begin{cases}
c_{30}(0) = 0.098178 & \mathrm{si\ } t=0, \\
c_{30}(t) \ \hbox{d\'efini en}\ \myref{c30} & \mathrm{si\ } t > t_1,
\end{cases}
\end{gather*}
o\`u $t_1$ est la plus petite partie imaginaire des z\'eros de
zeta~:~$\displaystyle{t_1=14.134725146}$.
\end{lem}
\begin{proof} \ \newline\noindent
Tout d'abord, remarquons que la somme \`a \'etudier prend sa valeur
maximale en $t_0=1$. Notons $\Sigma(t)$ cette nouvelle somme et $N(u)$ le nombre de z\'eros
non triviaux de $\zeta$ de partie imaginaire dans $[0,u]$.
\noindent
D'apr\`es le th\'eor\`eme de Backlund (cf. \cite{B}), $N(u)$ v\'erifie~:
$$
\left|N(u)-\frac{u}{2\pi}\log\left(\frac{u}{2\pi e}\right)\right|
 \le 0.137 \log u + 0.443 \log \log u + 5.225\,,\  \left(u\ge t_1\right).
$$
Nous en d\'eduisons cette in\'egalit\'e plus pratique pour
les applications num\'eriques que nous voulons mener par la suite~:
\begin{gather}
\label{dNu}
N_2(u) \le N(u) \le N_1(u)\,,\  u\ge t_1\\
\hbox{avec}\ N_1(u) = \frac{u}{2\pi}\log\left(\frac{u}{2\pi e}\right)+
0.29992 \log u + 5.225 \\
\hbox{et}\ N_2(u) = \frac{u}{2\pi}\log\left(\frac{u}{2\pi e}\right) -
0.29992 \log u - 5.225.
\end{gather}
\begin{itemize}
\item[-] Dans le cas o\`u $t=0$, comme il n'y a pas de z\'eros de zeta de partie
imaginaire inf\'erieure \`a $\displaystyle{t_1}$, nous avons l'\'egalit\'e ~:
$$ \Sigma (0) = \sum_{|\gamma|\ge 1} \frac1{|\gamma|^2} =
\sum_{|\gamma|\ge t_1} \frac1{|\gamma|^2}.
$$
et les majorations successives suivantes ~:
\begin{equation}
\Sigma (0) \le  2 \int_{t_1}^{+\infty} \frac{d\,N(u)}{u^2}
= 4 \int_{t_1}^{+\infty} \frac{N(u)}{u^3}\,d\,u
\le
0.098178.
\end{equation}
\item[-] Dans le cas o\`u $t>0$, nous avons ~:
\begin{equation}
\label{t>0}
\Sigma (t) \le \int_{|u-t| \ge 1} \frac{d\,N\left(|u|\right)}{(u-t)^2}
= \int_{u \ge t+1} \frac{d\,N(|u|)}{(u-t)^2}\,+\,
\int_{u \le t-1} \frac{d\,N\left(|u|\right)}{(u-t)^2}.
\end{equation}
Nous rappelons que nous voulons effectuer les calculs pour
$t\ge T_0$. On a donc $t_1\le t-1$, ce qui annule la seconde int\'egrale pour les valeurs de $u$ dans
$[-t_1;t_1]$ et finalement nous estimons les trois
int\'egrales suivantes en utilisant \myref{dNu}~:
\begin{eqnarray*}
\int_{-\infty}^{-t_1}  \frac{d\,N(-u)}{(u-t)^2}
&=& 2\int_{t_1}^{+\infty}  \frac{N(u)}{(u+t)^3}\,d\,u
\le 2\int_{t_1}^{+\infty}  \frac{N_1(u)}{(u+t)^3}\,d\,u \\
\int_{t_1}^{t-1}  \frac{d\,N(u)}{(u-t)^2}
&\le& N\left(t-1\right) \,-\, 2 \int_{t_1}^{t-1}
\frac{N_2(u)}{(t-u)^3}\,d\,u \\
\int_{t+1}^{+\infty} \frac{d\,N(u)}{(u-t)^2}
&\le& - N\left(t+1\right) \,+\, 2 \int_{t+1}^{+\infty}\frac{N_1(u)}{(u-t)^3}\,d\,u.
\end{eqnarray*}
On note respectivement $\displaystyle{J_1(t)}$,
$\displaystyle{J_2(t)}$, $\displaystyle{J_3(t)}$ les trois termes de
droite ci-dessus et $\displaystyle{c_{30}(t)}$ la majoration suivante de leur somme ~:
\begin{equation}
\label{c30}
c_{30}(t) =
2\int_{t_1}^{+\infty}  \frac{N_1(u)}{(u+t)^3}\,d\,u \,-\,
2 \int_{t_1}^{t-1} \frac{N_2(u)}{(t-u)^3}\,d\,u \,+\,
2 \int_{t+1}^{+\infty} \frac{N_1(u)}{(u-t)^3}\,d\,u.
\end{equation}
Remarquons que $c_{30}(t)=\mathcal{O}\left(\log t\right)$.
En effet, en int\'egrant dans la relation \myref{c30} les approximations
$N_i(u) = \frac{u}{2\pi} \log u +\mathcal{O}\left(u\right)$,
$i=1,\,2$, nous obtenons~:
\begin{eqnarray*}
\int_{t_1}^{+\infty}  \frac{N_1(u)}{(u+t)^3}\,d\,u
&=& \frac{t \log\left(t+t_1\right)}{4\pi\left(t+t_1\right)^2}
\,+\, \mathcal{O}\left(\frac1t\right)
= \mathcal{O}\left(\log t\right) ,\\
\int_{t_1}^{t-1} \frac{N_2(u)}{(t-u)^3}\,d\,t
&=& \frac{t^3 \log\left(t-1\right)}{4\pi\left(t-t_1\right)^2}
\,+\, \mathcal{O}\left(\log t\right)
= \frac{t\log\left(t-1\right)}{4\pi} \,+\, \mathcal{O}\left(\log  t\right) ,\\
\int_{t+1}^{+\infty} \frac{N_1(u)}{(u-t)^3}\,d\,u
&=&  \frac{t\log\left(t+1\right)}{4\pi}
\,+\, \mathcal{O}\left(\log t\right)
\end{eqnarray*}
$$
\hbox{et}\
c_{30}(t) = \frac{t}{2\pi} \log\left(\frac{t+1}{t-1}\right)
\,+\, \mathcal{O}\left(\log t\right)
= \frac{1}{\pi} \, \frac{t}{t-1}
\,+\, \mathcal{O}\left(\log t\right) = \mathcal{O}\left(\log t\right).
$$
\end{itemize}
\end{proof}
\begin{prop}
\label{Prop2.7}
\begin{multline*}
\sum_{{\begin{substack}{\varrho \in Z(\zeta) \\
                       y_k \ge t_0}
      \end{substack}}}
\Big[ D(\sigma-\beta+iy_k) + D(\sigma-1+\beta+iy_k) \Big]
\\ \ge
-\Big( \frac{M(0)c_{30}(kT_0)}2  \eta^2 +\frac{(1+2\kappa)mc_{30}(kT_0)}{2\sigma_0-1}  \eta^3 \Big)
\end{multline*}
\end{prop}
\begin{proof} \ \newline\noindent
D'apr\`es la majoration de $\tilde F$ \'etablie au paragraphe
\ref{Section5.2}, nous avons~:
\begin{multline}
\label{3.11}
\sum_{{\begin{substack}{\varrho \in Z(\zeta) \\
                       y_k \ge t_0}
      \end{substack}}}
\Big[ D(\sigma-\beta+iy_k) + D(\sigma-1+\beta+iy_k) \Big]
\\ \ge g_1(\theta) \sum_{{\begin{substack}{\varrho \in Z(\zeta) \\
                       y_k \ge t_0}
      \end{substack}}}
 \Re \Big( \frac1{\sigma-\beta+iy_k} + \frac1{\sigma-1+\beta+iy_k}
\\ -  \ \frac{\kappa}{\sigma-\beta +\delta +iy_k}
-\frac{\kappa}{\sigma-1+\beta +\delta +iy_k} \Big)
 \\  - \sum_{{\begin{substack}{\varrho \in Z(\zeta) \\
                       y_k \ge t_0}
      \end{substack}}}
\Big( |H(\sigma-\beta, y_k)|  + |H(\sigma-1+\beta, y_k)|
\\+ \kappa\ |H(\sigma-\beta+\delta, y_k)| +\kappa\ |H(\sigma-1+\beta+\delta, y_k)|  \Big)
\end{multline}
Le terme g\'en\'eral de la premi\`ere somme est positif puisque
$\displaystyle{ \delta \ge \frac{\sqrt 5 -1}2 }$, d'apr\`es le lemme de Stechkin
(cf. lemme 2, \cite{Stechkin1}).
\noindent
Il reste \`a minorer la seconde somme de \myref{3.11}.
Or le lemme \ref{Lemme2.2} permet de majorer respectivement
$\displaystyle{|H(\sigma-\beta, y_k)|}$ par $\displaystyle{\frac{M(0)}{y_k^2}\,\eta^2}$, puis
$\displaystyle{|H(\sigma-1+\beta, y_k)|}$, $\displaystyle{|H(\sigma-\beta+\delta, y_k)|}$ et
$\displaystyle{|H(\sigma-1+\beta+\delta, y_k)|}$ par
$\displaystyle{\frac{m}{(\sigma_0-1/2)y_k^2}\,\eta^3}$.
L'in\'egalit\'e \myref{3.11} devient alors~ :
\begin{multline}
\label{3.12}
\sum_{{\begin{substack}{\varrho \in Z(\zeta) \\
                       y_k \ge t_0}
      \end{substack}}}
\Big[ D (\sigma-\beta+iy_k)+D (\sigma-1+\beta+iy_k)\Big]
\\ \ge - \Big[  M(0)\,\eta^2 + \frac{(1+2\kappa)m}{\sigma_0-1/2}\,\eta^3 \Big]
\sum_{{\begin{substack}{\varrho \in Z(\zeta) \\
                       y_k \ge t_0}
      \end{substack}}} \frac1{y_k^2}
\end{multline}
Le lemme \ref{Lemme4} nous fournit une majoration de la somme de
droite et donc \myref{3.12} devient l'in\'egalit\'e annonc\'ee.
\end{proof}
\noindent On finit la preuve de la proposition \ref{Prop1.5} en d\'eduisant
tout d'abord de la proposition \ref{Prop2.7} que
\begin{equation}
\label{3.14}
\sum_{k=0}^4 a_k \sum_{{\begin{substack}{\varrho \in Z(\zeta) \\
                       y_k \ge t_0}
      \end{substack}}}
\Big[ D(\sigma-\beta+iy_k) + D(\sigma-1+\beta+iy_k) \Big]  \ge - \mathcal{C}_{31}(\eta) - \mathcal{C}_{32}(\eta)\kappa
\end{equation}
\begin{equation}
\label{C31}
\hbox{avec }\ \mathcal{C}_{31}(\eta)
= \frac12 \Big[  M(0)\, \eta^2+ \frac{m}{\sigma_0-1/2}\, \eta^3 \Big] \sum_{k=0}^4
a_k c_{30}(kT_0)
\end{equation}
\begin{equation}
\label{C32}
\mathcal{C}_{32}(\eta) =
 \frac{m}{\sigma_0-1/2} \sum_{k=0}^4 a_k c_{30}(kT_0)\, \eta^3
\end{equation}
(Nous rappelons que $c_{30}$ est d\'efini au lemme \ref{Lemme4}.)\\\noindent
\noindent
La proposition \ref{Prop2.6}, les in\'egalit\'es \myref{3.10} et \myref{3.14} donnent finalement~:
$$
\sum_{k=0}^4 a_k \sum_{\varrho \in Z(\zeta)}
\Big[ D(\sigma-\beta+iy_k) + D(\sigma-1+\beta+iy_k) \Big]
 \ge a_1 \tilde F(\sigma-\beta_0,0) - \mathcal{C}_{3}(\eta)
$$
\begin{multline}
\label{C3}
\hbox{avec}\ \mathcal{C}_{3}(\eta) =
 a_1 \Bigg[ \Bigg( \Big(\frac1{\delta}+\frac1{1-\eta_0+\delta}\Big)g_1(\theta)\,\eta +
 \Big(\frac1{\delta^3}+\frac1{(1-\eta_0+\delta)^3}\Big)m\,\eta^3 \Bigg)\kappa
 \\- \Bigg( g_1(\theta)\,\eta-m\,\eta^3 \Bigg) \Bigg] + \mathcal{C}_{31}(\eta) + \mathcal{C}_{32}(\eta) \kappa
\end{multline}
Ainsi, $\mathcal{C}_{3}(\eta)$ s'\'ecrit sous forme polynômiale~:
\begin{multline}
\label{C3bis}
\qquad\qquad\qquad\qquad\qquad\qquad p_1 \eta + p_2 \eta^2 + p_3 \eta^3,\\
\hbox{o\`u}\quad
p_1=a_1
g_1(\theta)\left(\left(\frac1{\delta}+\frac1{1-\eta_0+\delta}\right)\kappa-1\right)\
,\quad
p_2 = \frac{M(0)}{2}\sum_{k=0}^4 a_k c_{30}(kT_0)
\\ \hbox{ et }\
p_3 = \frac{3m}{2\sigma_0-1}\sum_{k=0}^4 a_k c_{30}(kT_0) +
a_1 m \left(\left(\frac1{\delta^3}+\frac1{(1-\eta_0+\delta)^3}\right)\kappa-1 \right).
\end{multline}
Les conditions impos\'ees \`a $\kappa$ et $\delta$ en
\myref{hyp-kappa} impliquent que $p_1$ est n\'egative et $p_2$ et $p_3$ positives.\\ \noindent
Et avec les valeurs num\'eriques choisies au d\'ebut du
paragraphe, on a plus exactment~:
\begin{gather*}
p_1 = -54.957\, , \ p_2 = 344\,602.065\, , \ p_3 = 3\,384\,045.191.
\end{gather*}
\subsection{\'Etude de $\Delta _1$ - Preuve de la proposition
~\ref{Prop1.2}~: }
\label{Section2.3}\noindent
Rappelons que le terme \'etudi\'e est
$\Delta_1(s) = T_1(s) - \kappa\ T_1(s+\delta)$
, o\`u
$$T_1(s) = -\frac12\log \pi
+ \frac12 \Re \frac{\Gamma'}{\Gamma}\Big(\frac{s}2+1\Big)
$$
\begin{lem}
\label{delta1}
$$
\Delta_1(\sigma+ik\gamma_0) \le c_1(k)
$$
\begin{gather*}
\quad\hbox{avec}\ \
c_1(0) = -\frac{1-\kappa}2 \log \pi
+ \frac12 \frac{\Gamma'}{\Gamma}\Big(\frac32\Big) - \frac{\kappa}2
\frac{\Gamma'}{\Gamma}\Big(\frac{\sigma_0+\delta}2+1\Big)\,,
\\
c_1(k) = -\frac{1-\kappa}2 \log\frac{2\pi}{k} + \frac12 \min\left(r_2(\sigma_0+2,3,kT_0),r_3(\sigma_0+2,3,kT_0)\right) \ \hbox{ si }\ k\ge1.
\end{gather*}
\end{lem}
\noindent Nous rappelons que $r_2$ et $r_3$ ont \'et\'e d\'efinis en \myref{r2} et
\myref{r3}.
\begin{proof}\ \\ \noindent
\begin{enumerate}
\item[$\bullet$]
Si $k=0$, nous avons imm\'ediatement grâce \`a la croissance de $\Re
\frac{\Gamma'}{\Gamma}$ sur $[0;+\infty[$~:
$$
\Delta_1(\sigma)  \le -\frac{1-\kappa}2 \log \pi
+ \frac12 \frac{\Gamma'}{\Gamma}\Big(\frac{3}2\Big) - \frac{\kappa}2
\frac{\Gamma'}{\Gamma}\Big(\frac{\sigma_0+\delta}2+1\Big) = c_1(0)
$$
\item[$\bullet$]
Si $k\ge1$, utilisons le lemme \ref{Gamma1} en prenant
$x_0=\sigma_0+2$, $x_1=3$ et $y_0=kT_0$~:
\begin{multline*}
\Delta_1(\sigma+ik\gamma_0) \le
 \frac{1-\kappa}2 \log \gamma_0
+ \frac{1-\kappa}2 \log\frac{k}{2\pi}
\\+ \frac12\min(r_2(\sigma_0+2,3,kT_0),r_3(\sigma_0+2,3,kT_0))
 \le \frac{1-\kappa}2 \log \gamma_0 + c_1(k)
\end{multline*}
o\`u les $c_1(k)$ sont des constantes n\'egatives.
\end{enumerate}
\end{proof}
\noindent Ainsi, en sommant le lemme \ref{delta1} pour $k=0,1,2,3,4$~:
\begin{equation}
\label{c1}
\sum_{k=0}^4a_k \Delta _{1}(\sigma + ik\gamma_0) \le \frac{A(1-\kappa)}2 \log \gamma_0 + c_1
\quad\hbox{avec}\quad
c_1 = \sum_{k=0}^4 a_k c_1(k)
\end{equation}
et donc
$$
f(0)\ \sum_{k=0}^4a_k \Delta _{1}(\sigma + ik\gamma_0)
 \le \frac{A}2(1-\kappa)g_1(\theta)\ \eta \log \gamma_0 + \mathcal{C}_1(\eta)
$$
o\`u $\mathcal{C}_1$ est la fonction n\'egative donn\'ee par~:
\begin{equation}
\label{C1}
\mathcal{C}_1(\eta) = c_1 g_1(\theta)\ \eta \le -2718.913 \eta
\end{equation}
\subsection{\'Etude de $D(s-1)$ - Preuve de la proposition~\ref{Prop1.3}~:}
\label{Section2.4}\noindent
Rappelons que
$D(\sigma-1+it) = \tilde F(\sigma-1, t)- \kappa\ \tilde F(\sigma-1+ \delta, t)$.
\noindent
Nous montrons tout d'abord une proposition interm\'ediaire~:
\begin{prop}
\label{Prop2.5}
\begin{eqnarray}
D(\sigma-1+it) \le
\begin{cases}
\tilde F(\sigma-1,0) - (238.212\eta-5\,533.813\eta^3) & \mathrm{si}\ t=0 \\
( - 20.991 \eta + 1\,403.284\eta^2)/t^2 & \mathrm{si}\ t\ge T_0
\end{cases}
\end{eqnarray}
\end{prop}
\begin{proof} \ \newline\noindent
Nous utilisons simplement les majorations du lemme \ref{Lemme2.2}.
\noindent
Dans le cas o\`u $t=0$, nous avons alors ~:
\begin{multline}
\label{2.9}
D(\sigma-1) = \tilde F(\sigma-1,0)- \kappa\ \tilde F(\sigma-1+\delta ,0)  \\
 \le \tilde F(\sigma-1,0)-\ \kappa
\Big( \frac{g_1(\theta)}{\delta} \eta - \frac{m}{\delta^3}\eta^3 \Big)
\end{multline}
et dans le cas o\`u $t \ge T_0$~:
\begin{equation}
\label{2.9bis}
\tilde F(\sigma-1, t)- \kappa\ \tilde F( \sigma-1+\delta, t)
\le M(-1) \frac{\eta^2}{t^2}
- \Big(g_1(\theta)\frac{\sigma_0-1+\delta}2\eta-m\eta^2\Big)\frac{\kappa}{t^2}
\end{equation}
Par cons\'equent, avec les valeurs choisies en d\'ebut de paragraphe,
nous avons la majoration explicite ~:~
$D(\sigma-1+it) \le (- 20.991 \eta +1403.284 \eta^2)/t^2 $,
ce qui termine la preuve de la proposition.
\end{proof}
\noindent Finalement, en prenant $t=k\gamma_0$ dans \myref{2.9} et \myref{2.9bis} avec successivement $k=0,1,2,3,4$,
la proposition \ref{Prop2.5} permet d'achever la preuve de la
proposition~\ref{Prop1.3}~:
$$
\sum_{k=0}^4 a_k D(\sigma-1+ik\gamma_0) \le a_0 \tilde F(\sigma-1,0) +
\mathcal{C}_2(\eta)
$$
o\`u $\mathcal{C}_2(\eta)$ s'\'ecrit sous la forme polyn\^omiale
\begin{multline}
\label{C2}
\qquad\qquad\qquad\qquad\qquad\qquad\mathcal{C}_2(\eta) = q_1\,\eta + q_2\,\eta^2+q_3\,\eta^3\,,\\
\hbox{avec}\
q_1 =   -\kappa \left( a_0 \frac{g_1(\theta)}{\delta} + \frac{\delta
  g_1(\theta)}2 \sum_{k=1}^4 \frac{a_k}{(kT_0)^2} \right) \le 0\,,\\
q_2 = \left( M(-1)+\kappa m \right)  \sum_{k=1}^4
\frac{a_k}{(kT_0)^2} \ge 0
\ \hbox{ et }\
q_3 =  a_0 \frac{m}{\delta^3}\kappa \ge 0\,,\qquad\qquad
\end{multline}
\begin{gather*}
\hbox{ et ici }\
q_1 = -1\,141.389 \, , \
q_2 = 2.794 \cdot10^{-15}\, , \
q_3 = 26\,515.117.
\end{gather*}
%
\subsection{\'Etude du reste $\Delta _2$
- Preuve de la proposition \ref{Prop1.4}}
\label{Section2.6}\noindent
Rappelons que $\Delta_2(s) = T_2(s)-\kappa\ T_2(s+\delta)$, avec~·
$$
T_2(s) = \frac{1}{2\pi}\int_{-\infty}^{+\infty}
\Re\frac{\Gamma'}{\Gamma}\Big(\frac14 + i\frac T2 \Big)
H(\sigma-1/2,t-T)d\,T + H(\sigma,t)
$$
\begin{lem}
\label{LemmeD2}
$$
\Delta_2(\sigma+ki\gamma_0) \le \mathcal{C}_4(\eta,k)
$$
$$\hbox{avec}\quad
\mathcal{C}_4(\eta,k)=\mathcal{C}_{41}(\eta,k)+\mathcal{C}_{42}(\eta,k)
$$
o\`u $\mathcal{C}_{41}(\eta,k)$ et $\mathcal{C}_{42}(\eta,k)$ sont d\'efinis en
\myref{C41} et \myref{C42}.
\end{lem}
\begin{proof} \ \newline\noindent
\begin{enumerate}
\item
\'Etudions d'abord le terme int\'egral.
\begin{multline*}
\Big| \int_{-\infty}^{+\infty}
\Re\frac{\Gamma'}{\Gamma}\Big(\frac14 +i\frac{T}2 \Big)
\Re \frac{F_2(x-i(T-y))}{(x-i(T-y))^2} d\,T \Big| \\
\le \int_{-\infty}^{+\infty}
\Big| \Re\frac{\Gamma'}{\Gamma}\Big(\frac14 +i\frac{T}2 \Big) \Big|
\Big| \Re \frac{F_2(x-i(T-y))}{(x-i(T-y))^2} \Big| d\,T
\end{multline*}
Comme
\begin{multline*}
\Big| \Re \frac{F_2(\sigma-1/2-i(T-t))}{(\sigma-1/2-i(T-t))^2} \Big|
= \eta^2 \Big| \int_{0}^{d_1(\theta)} \Re \frac{h"(t) e^{-(x-i(T-y))t/\eta}}{(x-i(T-y))^2} d\,t \Big|
\\ = \eta^2 \Big| \int_{0}^{d_1(\theta)} h"(t)\frac{e^{-xt/\eta}}{x^2+(T-y)^2} d\,t
\Big|
 \le \eta^2 \int_{0}^{d_1(\theta)} |h"(t)| \frac{e^{-xt/\eta}}{x^2+(T-y)^2} d\,t
\end{multline*}
D'apr\`es le th\'eor\`eme de Fubini, nous avons~:
\begin{multline*}
\int_{-\infty}^{+\infty}
\Big|\Re\frac{\Gamma'}{\Gamma}\Big(\frac14 +i\frac{T}2 \Big) \Big|
\Big|\Re \frac{F_2(x-i(T-y))}{(x-i(T-y))^2} \Big| d\,T  \\
\le \eta^2 \int_{-\infty}^{+\infty}
\Big| \Re\frac{\Gamma'}{\Gamma}\Big(\frac14 +i\frac{T}2 \Big) \Big|
\int_{0}^{d_1(\theta)} \frac{|h"(t)| e^{-xt/\eta}}{x^2+(T-y)^2} d\,t d\,T \\
= \eta^2 M\Big(\frac{x}{\eta}\Big) \int_{-\infty}^{+\infty}
\Big| \Re\frac{\Gamma'}{\Gamma}\Big(\frac14 +i\frac{T}2 \Big)
\Big| \frac1{x^2+(T-y)^2} d\,T \\
\end{multline*}
Enfin, en majorant $\Big|\Re\frac{\Gamma'}{\Gamma}\Big|$ grâce au lemme
\ref{Gamma2} et
$M\big(\frac{x}{\eta}\big)$ par $\frac{m}{x}\eta$ nous obtenons~:
\begin{multline*}
\label{C40}
\frac1{2\pi} \Big| \int_{-\infty}^{+\infty}
\Re\frac{\Gamma'}{\Gamma}\Big(\frac14 +i\frac{T}2 \Big)
\Re \frac{F_2(x-i(T-y))}{(x-i(T-y))^2} d\,T \Big| \\
\le \frac{m\eta^3}{2\pi x} \int_{-\infty}^{+\infty}
\frac{U_0(T)}{x^2+(T-y)^2} d\,T
= C_{40}(\eta,x,y)
\end{multline*}
Notons que si l'int\'egrale qui intervient est de l'ordre de
$\log t$, le $\eta^3$ qui la pr\'ec\`ede est de l'ordre de
$1/\log^3t$ ce qui fait qu'il est facile de montrer que cette
quantit\'e est d\'ecroissante en $t$ et qu'elle est donc major\'ee par $C_{40}(\eta,x,kT_0)$.
Notons
\begin{equation}
\label{C41}
C_{41}(\eta,k) = C_{40}(\eta,\sigma_0-1/2,kt_1)+\kappa\
C_{40}(\eta,\sigma_0-1/2+\delta,kt_1)
\end{equation}
\item
Il nous reste \`a approcher $H(\sigma,k\gamma_0)$ grâce au lemme \ref{Lemme2.2}~:
\begin{equation}
\label{c45k}
|H(\sigma,k\gamma_0)| + \kappa\
|H(\sigma+\delta,k\gamma_0)|
\le
C_{42}(\eta,k)
\end{equation}
\begin{equation}
\label{C42}
\hbox{avec}\quad
C_{42}(\eta,k) =
\begin{cases}
\Big(\frac1{\sigma_0^3}+\frac{\kappa}{(\sigma_0+\delta)^3}\Big) m\eta^3
&\ \mathrm { si }\ k=0 \\
\Big(\frac1{\sigma_0}+\frac{\kappa}{\sigma_0+\delta}\Big)\frac{m\eta^3}{(kT_0)^2}
& \mathrm { sinon }
\end{cases}
\end{equation}
\end{enumerate}
Finalement~:
$$
\sum_{k=0}^4 a_k \Delta_2(\sigma+ik\gamma_0) \le \mathcal{C}_4(\eta)
$$
\begin{equation}
\label{C4}
\hbox{ o\`u }\
\mathcal{C}_4(\eta) =  \sum_{k=0}^4 a_k
\Big(\mathcal{C}_{41}(\eta,k)+\mathcal{C}_{42}(\eta,k)\Big)\
\hbox{ est toujours positive.}
\end{equation}
$$
\mathcal{C}_4(\eta) \le 2.3887 \cdot10^6 \eta^3.
$$
\end{proof}
\ \\ \\ \noindent
Nous d\'etaillons ci-dessous les \'etapes cons\'ecutives des
calculs.
Nous donnons les valeurs successives prises pour $r$ et
$R$, ainsi que celles des
param\`etres $\eta_0$, $\kappa$ et $\delta$ impliqu\'es et enfin celles trouv\'ees pour
$R_0$. \\\noindent
Nous donnons aussi le terme reste
$$\mathcal{C}(\eta) = \mathcal{C}_1(\eta)
+ \mathcal{C}_2(\eta) + \mathcal{C}_3(\eta) + \mathcal{C}_4(\eta) = \alpha_1\, \eta +
\alpha_2 \,\eta^2 + \alpha_3\, \eta^3\,, $$
o\`u $\alpha_1$ est
n\'egative sous la condition \myref{hyp-kappa}, et $\alpha_2$ et $\alpha_3$ sont toujours positives. Ainsi, comme on l'a
expliqu\'e au paragraphe \ref{Section4}, $\mathcal{C}(\eta)$ est n\'egatif sur
$[0;\eta_0]$ lorsque $\mathcal{C}(\eta_0)$ est n\'egatif, et il n'influe donc
pas sur la valeur finale de $R_0$.
$$
\begin{array}{|c|c|c|c|c|c|}
\hline
Etape & R & r & \eta_0 \cdot{10^3} & \kappa & \delta \\
\hline
1 & 9.645908801 & 5.97484 & 7.63319 & 0.438904 & 0.620626 \\ \hline
2 & 5.974849075 & 5.73045 & 7.95873 & 0.438525 & 0.620748 \\ \hline
3 & 5.730454010 & 5.70487 & 7.99441 & 0.438483 & 0.620762 \\ \hline
4 & 5.704872616 & 5.70208 & 7.99832 & 0.438479 & 0.620763 \\ \hline
5 & 5.702089881 & 5.70178 & 7.99874 & 0.438478 & 0.620763 \\ \hline
6 & 5.701785245 & 5.70174 & 7.99880 & 0.438478 & 0.620763 \\ \hline
\end{array}
$$
$$
\begin{array}{|c|c|c|c|c|c|}
\hline
Etape
& \alpha_1 & \alpha_2 & \alpha_3 & \mathcal{C}(\eta_0)  & R_0 \\ \hline
1 & -3\,915.260 & 344\,602.065 & 5\,799\,250.773 & -7.22827& 5.974849075 \\ \hline
2 & -3\,916.747 & 344\,602.065 & 5\,841\,345.585 &-7.22089 & 5.730454010 \\ \hline
3 & -3\,916.907 & 344\,602.065 & 5\,846\,103.683 & -6.30271& 5.704872616\\ \hline
4 & -3\,916.907 & 344\,602.065 & 5\,846\,103.683 & -6.29209& 5.702089881 \\ \hline
5 & -3\,916.926 & 344\,602.065 & 5\,846\,682.864 & -6.29080& 5.701785245\\ \hline
6 & -3\,916.927 & 344\,602.065 & 5\,846\,689.069 &-6.29065 & 5.701752890 \\ \hline
\end{array}
$$

\ \\\\

\footnotesize

\parbox{3.2in}{

\noindent Habiba Kadiri\\
{\sc D\'epartement de Math\'ematiques et Statistique\\
Universit\'e de Montr\'eal}\\
CP 6128 succ Centre-Ville\\
Montr\'eal QC  H3C 3J7\\
Qu\'ebec, Canada\\
e-mail: kadiri@dms.umontreal.ca}

\end{document}